\documentclass[12pt]{article}
\usepackage{e-jc}
\newtheorem{Theorem}{Theorem}
\newtheorem{Lemma}[Theorem]{Lemma}

\newtheorem{Proposition}[Theorem]{Proposition}
\newtheorem{Definition}[Theorem]{Definition}
\date{\dateline{June 2, 2009}{Oct 13, 2009}}
\begin{document}
\title{The evolution of uniform random planar graphs}
\author{Chris Dowden}
\maketitle

\begin{abstract}

Let $P_{n,m}$ denote the graph
taken uniformly at random from the set of all planar graphs on $\{1,2, \ldots, n \}$
with exactly $m(n)$ edges.
We use counting arguments to investigate the probability that $P_{n,m}$
will contain given components and subgraphs,
finding that there is different asymptotic behaviour depending on the ratio $\frac{m}{n}$.

\end{abstract}

\section{Introduction} \label{aintro}

Random planar graphs have recently been the subject of much activity,
and many properties of the standard random planar graph $P_{n}$
(taken uniformly at random from the set of all planar graphs on $\{1,2, \ldots, n\}$)
are now known.
For example,
in~\cite{mcd} it was shown that $P_{n}$ will asymptotically almost surely
(a.a.s., that is, with probability tending to $1$ as $n$ tends to infinity)
contain at least linearly many copies of any given planar graph.
By combining the counting methods of~\cite{mcd}
with some rather precise results of~\cite{gim}, obtained using generating functions,
the exact limiting probability for the event that $P_{n}$
will contain any given component is also known.

More recently,
attention has turned to the graph $P_{n,m}$
taken uniformly at random from the set $\mathcal{P}(n,m)$ of all planar graphs on $\{1,2,\ldots,n\}$
with exactly $m(n)$ edges.
It is well known that we must have $m<3n$ for planarity to be possible
and also that $P_{n,m}$ behaves in exactly the same way as the general random graph $G_{n,m}$
if $\frac{m}{n} < \frac{n}{2} - \omega \left( n^{2/3} \right)$,
so the interest lies with the region $\frac{1}{2} \leq \liminf \frac{m}{n} \leq \limsup \frac{m}{n} \leq 3$.

In~\cite{ger},
the case when $m = \lfloor qn \rfloor$ for fixed $q$ was investigated
and it was shown that, for all $q \in (1,3)$,
$P_{n, \lfloor qn \rfloor}$ will a.a.s.~contain at least linearly many copies of any given planar graph,
as with $P_{n}$.
It was also shown that
the probability that $P_{n, \lfloor qn \rfloor}$ will contain an isolated vertex is bounded away from $0$ as $n \to \infty$
(for all $q$)
and hence that
the probability that $P_{n, \lfloor qn \rfloor}$ will be connected is bounded away from $1$.
In~\cite{gim},
the precise limiting probability for
$\mathbf{P}[P_{n, \lfloor qn \rfloor}$ will be connected]
was then determined,
using generating functions.

In this paper,
we use counting arguments to extend the current knowledge of $P_{n,m}$.
We investigate the probability that $P_{n,m}$ will contain given components
and the probability that $P_{n,m}$ will contain given subgraphs,
both for general $m(n)$,
and show that there is different behaviour
depending on which `region' the ratio $\frac{m(n)}{n}$ falls into.
Hence, this change as $\frac{m}{n}$ varies can be thought of as the `evolution'
of uniform random planar graphs. \\

We will start in Section~\ref{pen}
by collecting up various lemmas on $P_{n,m}$
that shall prove useful to us.
In Section~\ref{cptlow},
we will then obtain lower bounds for
$\mathbf{P} := \mathbf{P}[P_{n,m}$ will have a component isomorphic to $H]$
(where by `lower bound' we mean a result such as
$\liminf \mathbf{P} > 0$ or $\mathbf{P} \to 1$),
and in Section~\ref{cyccpt} we will obtain exactly complementary upper bounds.
Finally, in Section~\ref{sub},
we will then look at
the probability that
$P_{n,m}$ will have a copy of $H$ (i.e.~any subgraph isomorphic to $H$).

A summary of our results is given in Tables~\ref{cpttab} and~\ref{subtab}.
These tables both have three columns,
corresponding to the sign of $e(H) - |H|$
(the excess of edges over vertices),
and also different rows,
corresponding to the size of $\frac{m(n)}{n}$.
We use $\underline{\lim}$ to denote $\liminf$
and $\overline{\lim}$ to denote $\limsup$,
and `T\ref{gen3}' (for example) refers to Theorem~\ref{gen3}.

\begin{table} [ht]
\begin{tabular}{|c|c|c|c|}
\hline
&
$e(H)<|H|$ &
$e(H)=|H|$ &
$e(H) > |H|$ \\
\hline
$0 < \textrm{\small{$\underline{\lim}$} } \frac{m}{n}$ &
$\mathbf{P} \to 1$ (Thm \ref{tree6}) &
\small{$\underline{\lim}$} $\mathbf{P} > 0$ (T\ref{gen4}) &
$\mathbf{P} \to 0$ (Thm \ref{cyc41}) \\
\& $\frac{m}{n} \leq 1+o(1)$ &
&
\small{$\overline{\lim}$} $\mathbf{P} < 1$ (T\ref{cyc52}) &
\\
\hline
$1 < \textrm{\small{$\underline{\lim}$} } \frac{m}{n}$ &
\small{$\underline{\lim}$} $\mathbf{P} > 0$ (T\ref{gen3}) &
\small{$\underline{\lim}$} $\mathbf{P} > 0$ (T\ref{gen3}) &
\small{$\underline{\lim}$} $\mathbf{P} > 0$ (T\ref{gen3}) \\
\& $\textrm{\small{$\overline{\lim}$} } \frac{m}{n} < 3$ &
\small{$\overline{\lim}$} $\mathbf{P} < 1$ (T\ref{conn4}) &
\small{$\overline{\lim}$} $\mathbf{P} < 1$ (T\ref{conn4}) &
\small{$\overline{\lim}$} $\mathbf{P} < 1$ (T\ref{conn4}) \\
\hline
$\frac{m}{n} \to 3$ &
$\mathbf{P} \to 0$ (Thm \ref{conn5}) &
$\mathbf{P} \to 0$ (Thm \ref{conn5}) &
$\mathbf{P} \to 0$ (Thm \ref{conn5}) \\
\hline
\end{tabular}
\caption{A description of $\mathbf{P} := \mathbf{P}[P_{n,m}$ will have a component isomorphic to $H]$.} \label{cpttab}
\end{table}

\begin{table} [ht]
\begin{tabular}{|c|c|c|c|}
\hline
&
$e(H)<|H|$ &
$e(H)=|H|$ &
$e(H) > |H|$ \\
\hline
$0 < \underline{\lim}~\frac{m}{n}$ &
$\mathbf{P}^{\prime} \to 1$ (Thm \ref{tree6}) &
\small{$\underline{\lim}$} $\mathbf{P}^{\prime} > 0$ (T\ref{gen4}) &
$\mathbf{P}^{\prime} \to 0$ (Thm \ref{msub3}) \\
\& $\overline{\lim}~\frac{m}{n} < 1$ &
&
\small{$\overline{\lim}$} $\mathbf{P}^{\prime} < 1$ (T\ref{cyc54}) &
\\
\hline
$\frac{m}{n} \to 1$ &
$\mathbf{P}^{\prime} \to 1$ (Thm \ref{tree6}) &
$\mathbf{P}^{\prime} \to 1$ (Thm \ref{unisub1}) &
Unknown  \\
\hline
$\underline{\lim}~\frac{m}{n} > 1$ &
$\mathbf{P}^{\prime} \to 1$ (Thm \ref{sub4}) &
$\mathbf{P}^{\prime} \to 1$ (Thm \ref{sub4}) &
$\mathbf{P}^{\prime} \to 1$ (Thm \ref{sub4}) \\
\hline
\end{tabular}
\caption{A description of $\mathbf{P}^{\prime} :=
\mathbf{P}[P_{n,m}$ will have a copy of $H]$.} \label{subtab}
\end{table}

\section{Appearances, Pendant Edges \& Addable Edges} \label{pen}

In this section,
we shall lay the groundwork for the rest of the paper by noting some useful properties of $P_{n,m}$.
We will see results on the number of `appearances' (special subgraphs) in $P_{n,m}$,
the number of `pendant' edges
(i.e.~edges incident to a vertex of degree $1$),
and the number of `addable' edges
(i.e.~edges that can be added to $P_{n,m}$ without violating planarity).
All of these shall be important ingredients in the counting arguments of later sections. \\

We start with the definition of an appearance:

\begin{Definition} \label{defapps}
Let $H$ be a graph on the vertex set $\{1,2,\ldots,|H|\},$
and let $G$ be a graph on the vertex set $\{1,2,\ldots,n\}$,
where $n>|H|$.
Let $W \subset V(G)$ with $|W|=|H|$,
and let the `root' $r_{W}$ denote the least element in $W$.
We say that $H$ \emph{appears} at $W$ in $G$ if
(a) the increasing bijection from ${1,2,\ldots,|H|}$ to $W$ gives
an isomorphism between $H$ and the induced subgraph $G[W]$ of $G$;
and (b) there is exactly one edge in $G$ between $W$ and the rest of $G$,
and this edge is incident with the root $r_{W}$
(see Figure~\ref{app}).
We let $f_{H}(G)$ denote the number of appearances of $H$ in $G$,
that is the number of sets $W \subset V(G)$ such that $H$ appears at $W$ in $G$.
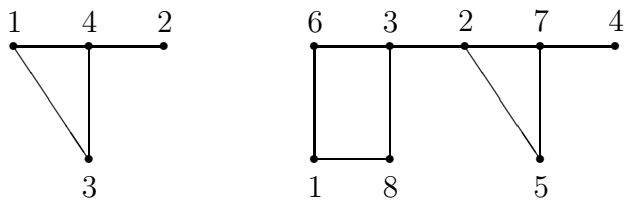
\begin{figure} [ht]
\setlength{\unitlength}{1cm}
\begin{picture}(10,2.5)(0,-0.5)
\put(2,1.5){\line(1,0){2}}
\put(1.91,1.7){$1$}
\put(2,1.5){\circle*{0.1}}
\put(2.91,1.7){$4$}
\put(3,1.5){\circle*{0.1}}
\put(3.91,1.7){$2$}
\put(4,1.5){\circle*{0.1}}
\put(3,1.5){\line(0,-1){1.5}}
\put(2,1.5){\line(2,-3){1}}
\put(2.91,-0.5){$3$}
\put(3,0){\circle*{0.1}}

\put(6,1.5){\line(1,0){4}}
\put(6,0){\line(1,0){1}}
\put(6,1.5){\line(0,-1){1.5}}
\put(7,1.5){\line(0,-1){1.5}}
\put(5.91,1.7){$6$}
\put(6,1.5){\circle*{0.1}}
\put(6.91,1.7){$3$}
\put(7,1.5){\circle*{0.1}}
\put(5.91,-0.5){$1$}
\put(6,0){\circle*{0.1}}
\put(6.91,-0.5){$8$}
\put(7,0){\circle*{0.1}}

\put(7.91,1.7){$2$}
\put(8,1.5){\circle*{0.1}}
\put(8.91,1.7){$7$}
\put(9,1.5){\circle*{0.1}}
\put(9.91,1.7){$4$}
\put(10,1.5){\circle*{0.1}}
\put(9,1.5){\line(0,-1){1.5}}
\put(8,1.5){\line(2,-3){1}}
\put(8.91,-0.5){$5$}
\put(9,0){\circle*{0.1}}

\end{picture}
\caption{A graph $H$ and an appearance of $H$.} \label{app}
\end{figure}
\end{Definition}

The following result on appearances was given in~\cite{ger}:

\begin{Proposition}[\cite{ger}, Theorem 3.1] \label{ger T3.1}
Let $H$ be a (fixed) connected planar graph on the vertices $\{1,2, \ldots ,|H|\}$
and let $q \in (1,3)$ be a constant.
Then there exists a constant $\alpha(H,q)>0$ such that
\begin{displaymath}
\mathbf{P} \left[f_{H} \left( P_{n, \lfloor qn \rfloor} \right) \leq \alpha n \right] = e^{- \Omega (n)}.
\end{displaymath}
\end{Proposition}

It is, in fact,
fairly easy to deduce from the proof of Proposition~\ref{ger T3.1} given in~\cite{ger}
that the result actually holds uniformly in $q$
(see~\cite{dow} for details).
Hence,
we may actually obtain the following stronger version:

\begin{Lemma}[\cite{dow}, Lemma 13] \label{gen2}
Let $H$ be a (fixed) connected planar graph on the vertices $\{1,2, \ldots ,|H|\}$,
let $b>1$ and $B<3$ be constants,
and let $m(n) \in [bn,Bn]$ for all $n$.
Then there exists a constant $\alpha = \alpha(H,b,B) >0$ such that
\begin{displaymath}
\mathbf{P}[f_{H}(P_{n,m}) \leq \alpha n] = e^{- \Omega (n)}.
\end{displaymath}
\end{Lemma}

An important consequence of Lemma~\ref{gen2} is that
$P_{n,m}$ will a.a.s.~contain a copy of any given planar graph if
$1< \liminf \frac{m}{n} \leq \limsup \frac{m}{n} < 3$
(as noted in the introduction).
However,
it is the precise uncomplicated structure of appearances themselves
that will be particularly useful to us during this paper. \\

It follows from Lemma~\ref{gen2}
(with $H$ as an isolated vertex)
that $P_{n,m}$ will a.a.s.~have linearly many pendant edges if
$1 < \liminf \frac{m}{n} \leq \limsup \frac{m}{n} < 3$.
It is fairly intuitive that we ought to be able to drop the lower bound
to $\liminf \frac{m}{n} > 0$ for this particular case,
and this is indeed shown in~\cite{dow}.
Hence, we obtain:

\begin{Lemma}[\cite{dow}, Theorem 16] \label{pen5}
Let $b>0$ and $B<3$ be constants and
let $m(n) \in [bn,Bn]$ for all $n$.
Then there exists a constant $\alpha = \alpha(b,B) > 0$ such that
\begin{displaymath}
\mathbf{P}[P_{n,m} \textrm{ will have less than } \alpha n \textrm{ pendant edges }] = e^{- \Omega (n)}.
\end{displaymath}
\end{Lemma}

\phantom{p}

We now move on to the final topic of this section, that of `addable' edges:

\begin{Definition}
Given a planar graph $G$, we call a non-edge $e$ \emph{addable} in $G$
if the graph $G+e$ obtained by adding $e$ as an edge is still planar.
We let \emph{add}$(G)$ denote the set of addable non-edges of $G$
(note that the graph obtained by adding two edges in \emph{add}$(G)$ may well not be planar)
and we let \emph{add}$(n,m)$
denote the minimum value of $|\emph{add}(G)|$ over all graphs $G \in \mathcal{P}(n,m)$.
\end{Definition}

In future sections,
we shall often wish to choose an edge to insert into a graph without violating planarity,
and we will want to know how many choices we have.
A very helpful result is given implicitly in Theorem 1.2 of~\cite{ger2}:

\begin{Lemma}[\cite{ger2}, Theorem 1.2] \label{add41}
Let $m(n) \leq (1+o(1))n$.
Then
\begin{displaymath}
\emph{add}(n,m)=\omega(n),
\end{displaymath}
i.e.~$\frac{\textrm{\small{$\emph{add}$}}(n,m)}{n} \to \infty$ as $n \to \infty$.
\end{Lemma}

We should also note that a useful higher estimate for the case
$\limsup \frac{m}{n} < 1$ can be obtained very simply:

\begin{Lemma} \label{add2}
Let $A < 1$ be a constant and let $m(n) \leq An$ for all $n$.
Then
\begin{displaymath}
\emph{add}(n,m) \geq (1+o(1)) \left( \frac{(1-A)^{2}}{2} \right) n^{2}.
\end{displaymath}
\end{Lemma}
\textbf{Proof}
Any graph in $\mathcal{P}(n,m)$ must have at least $n-m=(1-A)n$ components,
and it is known that inserting an edge between any two vertices in different components will not violate planarity.
Hence, add$(n,m) \geq \left( ^{(1-A)n} _{\phantom{ww}2} \right)$.
\phantom{www}
\setlength{\unitlength}{.25cm}
\begin{picture}(1,1)
\put(0,0){\line(1,0){1}}
\put(0,0){\line(0,1){1}}
\put(1,1){\line(-1,0){1}}
\put(1,1){\line(0,-1){1}}
\end{picture} \\

\section{Components I: Lower Bounds} \label{cptlow}

We now come to the first main section of this paper,
where we shall start to use the results of Section~\ref{pen} to investigate
$\mathbf{P} := \mathbf{P}[P_{n,m}$
will have a component isomorphic to $H]$.
We shall first see (in Theorem~\ref{gen3})
that $\liminf \mathbf{P} > 0$ for all connected planar $H$
if $1 < \liminf \frac{m}{n} \leq \limsup \frac{m}{n} < 3$,
then (in Theorem~\ref{gen4})
that the lower bound on $\frac{m}{n}$ can be reduced to $\liminf \frac{m}{n} > 0$
if $e(H) \leq |H|$,
and thirdly (in Theorem~\ref{tree6}) that
$\mathbf{P} \to 1$ if
$0 < \liminf \frac{m}{n} \leq \limsup \frac{m}{n} \leq 1$
and $H$ is a tree.
Finally, we will show (in Theorem~\ref{tree4})
that $P_{n,m}$ will a.a.s.~have \textit{linearly} many components isomorphic to any given tree if
$0 < \liminf \frac{m}{n} \leq \limsup \frac{m}{n} < 1$. \\

We start with our aforementioned result for general connected planar $H$:

\begin{Theorem} \label{gen3}
Let $H$
be a (fixed) connected planar graph,
let $b>1$ and $B<3$ be constants, and let $m(n) \in [bn,Bn]$ for all $n$.
Then there exist constants $\epsilon (H,b,B) > 0$
and $N(H,b,B)$ such that
\begin{displaymath}
\mathbf{P}[\textrm{$P_{n,m}$
will have a component isomorphic to $H$}  ]
\geq \epsilon \textrm{ for all } n \geq N.
\end{displaymath}
\end{Theorem}
\textbf{Sketch of Proof}
We shall suppose that the result is false.
Thus,
there exist arbitrarily large values of $n$ for which a typical graph in $\mathcal{P}(n,m)$
will have no components isomorphic to $H$,
but will have many appearances of $K_{4}$
(by Lemma~\ref{gen2}).
From each such graph,
we shall construct graphs in $\mathcal{P}(n,m)$
that do have a component isomorphic to $H$.

We start by deleting edges from some of our appearances of $K_{4}$ to create isolated vertices,
on which we then build a component isomorphic to $H$.
By inserting extra edges in appropriate places elsewhere,
we hence obtain graphs that are also in $\mathcal{P}(n,m)$.
The fact that the original graphs contained
no components isomorphic to $H$
can then be used to show that there isn't too much double-counting,
and so we find that we have actually constructed a decent number of \textit{distinct} graphs in $\mathcal{P}(n,m)$
that have components isomorphic to $H$,
which is what we wanted to prove.
\\
\\
\textbf{Full Proof}
Let $\epsilon \in (0,1)$.
Since $\frac{m}{n} \in [b,B]$ for all $n$,
by Lemma~\ref{gen2}
there exist constants $\alpha = \alpha(b,B) >0$ and $N(b,B)$ such that, for all $n \geq N$,
$\mathbf{P}[P_{n,m}$ will have at least $\alpha n$ appearances of $K_{4}]
\geq \left( 1- \frac{\epsilon}{2} \right)$.
Note that any appearances of $K_{4}$ must be vertex-disjoint,
by $2$-edge-connectedness.

Consider an $n \geq N$ and suppose that
$\mathbf{P}[P_{n,m}$
will have a component isomorphic to
$H] < 1-\epsilon$
(if not, then we are certainly done).
Let $\mathcal{G}_{n}$ denote the set of graphs in $\mathcal{P}(n,m)$
with (i) no components isomorphic to $H$
and (ii) at least $\alpha n$ vertex-disjoint appearances of $K_{4}$.
Then, under our assumption,
we have $|\mathcal{G}_{n}| \geq \frac{\epsilon}{2} |\mathcal{P}(n,m)|$.
We shall use $\mathcal{G}_{n}$ to construct graphs in $\mathcal{P}(n,m)$
that do have a component isomorphic to $H$.

\begin{figure} [ht]
\setlength{\unitlength}{0.785cm}
\begin{picture}(20,3.3)(-0.25,-0.3)

\put(0,0){\line(1,0){0.9}}
\put(0,0){\line(3,2){0.45}}
\put(0,0){\line(3,5){0.45}}
\put(0.45,0.3){\line(0,1){0.45}}
\put(0.9,0){\line(-3,2){0.45}}
\put(0.9,0){\line(-3,5){0.45}}
\put(0,0){\circle*{0.1}}
\put(0.9,0){\circle*{0.1}}
\put(0.45,0.3){\circle*{0.1}}
\put(0.45,0.75){\circle*{0.1}}

\put(1.2,0){\line(1,0){0.9}}
\put(1.2,0){\line(3,2){0.45}}
\put(1.2,0){\line(3,5){0.45}}
\put(1.65,0.3){\line(0,1){0.45}}
\put(2.1,0){\line(-3,2){0.45}}
\put(2.1,0){\line(-3,5){0.45}}
\put(1.2,0){\circle*{0.1}}
\put(2.1,0){\circle*{0.1}}
\put(1.65,0.3){\circle*{0.1}}
\put(1.65,0.75){\circle*{0.1}}

\put(2.4,0){\line(1,0){0.9}}
\put(2.4,0){\line(3,2){0.45}}
\put(2.4,0){\line(3,5){0.45}}
\put(2.85,0.3){\line(0,1){0.45}}
\put(3.3,0){\line(-3,2){0.45}}
\put(3.3,0){\line(-3,5){0.45}}
\put(2.4,0){\circle*{0.1}}
\put(3.3,0){\circle*{0.1}}
\put(2.85,0.3){\circle*{0.1}}
\put(2.85,0.75){\circle*{0.1}}

\put(3.6,0){\line(1,0){0.9}}
\put(3.6,0){\line(3,2){0.45}}
\put(3.6,0){\line(3,5){0.45}}
\put(4.05,0.3){\line(0,1){0.45}}
\put(4.5,0){\line(-3,2){0.45}}
\put(4.5,0){\line(-3,5){0.45}}
\put(3.6,0){\circle*{0.1}}
\put(4.5,0){\circle*{0.1}}
\put(4.05,0.3){\circle*{0.1}}
\put(4.05,0.75){\circle*{0.1}}

\put(2.25,2.25){\oval(5,1.5)}

\put(0.45,0.75){\line(0,1){0.75}}
\put(0.45,1.5){\circle*{0.1}}
\put(1.65,0.75){\line(0,1){0.75}}
\put(1.65,1.5){\circle*{0.1}}
\put(2.85,0.75){\line(0,1){0.75}}
\put(2.85,1.5){\circle*{0.1}}
\put(4.05,0.75){\line(0,1){0.75}}
\put(4.05,1.5){\circle*{0.1}}

\put(6,1.5){\vector(1,0){1}}

\put(8.95,0.3){\line(0,1){0.45}}
\put(9.4,0){\line(-3,2){0.45}}
\put(9.4,0){\line(-3,5){0.45}}
\put(9.4,0){\circle*{0.1}}
\put(8.95,0.3){\circle*{0.1}}
\put(8.95,0.75){\circle*{0.1}}

\put(9.7,0){\line(1,0){0.9}}
\put(9.7,0){\line(3,5){0.45}}
\put(10.6,0){\line(-3,5){0.45}}
\put(9.7,0){\circle*{0.1}}
\put(10.6,0){\circle*{0.1}}
\put(10.15,0.75){\circle*{0.1}}

\put(10.9,0){\line(1,0){0.9}}
\put(10.9,0){\line(3,5){0.45}}
\put(11.8,0){\line(-3,5){0.45}}
\put(10.9,0){\circle*{0.1}}
\put(11.8,0){\circle*{0.1}}
\put(11.35,0.75){\circle*{0.1}}

\put(12.1,0){\line(3,2){0.45}}
\put(12.1,0){\line(3,5){0.45}}
\put(12.55,0.3){\line(0,1){0.45}}
\put(12.1,0){\circle*{0.1}}
\put(12.55,0.3){\circle*{0.1}}
\put(12.55,0.75){\circle*{0.1}}

\put(10.75,2.25){\oval(5,1.5)}

\put(8.95,0.75){\line(0,1){0.75}}
\put(8.95,1.5){\circle*{0.1}}
\put(10.15,0.75){\line(0,1){0.75}}
\put(10.15,1.5){\circle*{0.1}}
\put(11.35,0.75){\line(0,1){0.75}}
\put(11.35,1.5){\circle*{0.1}}
\put(12.55,0.75){\line(0,1){0.75}}
\put(12.55,1.5){\circle*{0.1}}

\put(14,0.75){\line(1,0){0.9}}
\put(14,0.75){\line(3,2){0.45}}
\put(14,0.75){\line(3,5){0.45}}
\put(14.45,1.05){\line(0,1){0.45}}
\put(14.9,0.75){\line(-3,2){0.45}}
\put(14,0.75){\circle*{0.1}}
\put(14.9,0.75){\circle*{0.1}}
\put(14.45,1.05){\circle*{0.1}}
\put(14.45,1.5){\circle*{0.1}}

\put(9.25,2){\line(0,1){0.5}}
\put(9.75,2){\line(0,1){0.5}}
\put(10.25,2){\line(0,1){0.5}}
\put(10.75,2){\line(0,1){0.5}}
\put(11.25,2){\line(0,1){0.5}}
\put(11.75,2){\line(0,1){0.5}}
\put(12.25,2){\line(0,1){0.5}}

\put(13.8,0.45){\small{$v_{3}$}}
\put(14.8,0.45){\small{$v_{4}$}}
\put(14.55,1.05){\small{$v_{2}$}}
\put(14.55,1.55){\small{$v_{1}$}}

\put(-0.2,-0.3){\small{$v_{1}$}}
\put(1.5,0.1){\tiny{$v_{2}$}}
\put(2.7,0.1){\tiny{$v_{3}$}}
\put(4.5,-0.3){\small{$v_{4}$}}

\put(14.25,0){$H$}

\end{picture}
\caption{Constructing a component isomorphic to $H$.} \label{cpt1}

\end{figure}
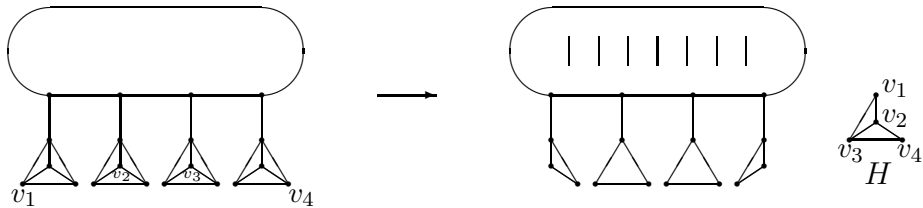

Consider a graph $G \in \mathcal{G}_{n}$.
We may assume that $n$ is large enough that $\alpha n \geq |H|$.
Thus, we may choose $|H|$ of the (vertex-disjoint) appearances of $K_{4}$ in $G$
$\left( \textrm{at least} \left( ^{\lceil \alpha n \rceil} _{\phantom{i}|H|} \right) \textrm{ choices} \right)$,
and for each of these chosen appearances we may choose a `special' vertex in the $K_{4}$ that is not the
root $\left( 3^{|H|} \textrm{ choices} \right)$.
Let us then delete all $3|H|$ edges that are incident to the `special' vertices and
insert edges between these $|H|$ newly isolated vertices in such a way that they now form a component
isomorphic to $H$ (see Figure~\ref{cpt1}).

To maintain the correct number of edges, we should insert $3|H|-e(H)$ extra ones somewhere into the graph,
making sure that we maintain planarity.
We will do this in such a way that we do not interfere with our new component
or with the chosen appearances of $K_{4}$
(which are now appearances of $K_{3}$).
Thus, the part of the graph where we wish to insert edges contains $n-4|H|$ vertices and $m - 7|H|$ edges.
We know that there exists a triangulation on these vertices containing these edges,
and clearly inserting an edge from this triangulation would not violate planarity.
Thus, we have at least
$\left( ^{3(n-4|H|)-6-(m -7|H|)} _{\phantom{wwww} 3|H|-e(H)} \right)$
choices for where to add the edges.

Therefore, in total we have at least
$ |\mathcal{G}_{n}|
\left( ^{\lceil \alpha n \rceil} _{\phantom{i}|H|} \right)
3^{|H|}
\left( ^{ 3(n-4|H|)-6-(m -7|H|) }_{\phantom{wwww} 3|H|-e(H) } \right)
= |\mathcal{G}_{n}|
\Theta \left( n^{4|H|-e(H)} \right)$
ways to build (not necessarily distinct) graphs in $\mathcal{P}(n,m)$
that have a component isomorphic to $H$.

We will now consider the amount of double-counting: \\
Each of our constructed graphs will contain at most $4|H|-e(H)+1$ components isomorphic to $H$
(since there were none originally;
we have deliberately built one; and
we may have created at most one extra one each time we cut a `special' vertex away from its $K_{4}$ or
added an edge in the rest of the graph).
Hence, we have at most $4|H|-e(H)+1$ possibilities
for which were our $|H|$ `special' vertices.
Since appearances of $K_{3}$ must be vertex-disjoint, by $2$-edge-connectedness,
we have at most $\frac{n}{3}$ of them and hence
at most
$\left( \frac{n}{3} \right)^{|H|}$
possibilities
for where the `special' vertices were originally.
There are then at most
$\left(^{m-e(H)-4|H|}_{\phantom{w}3|H|-e(H)}\right)$
possibilities for which edges were added
in the rest of the graph
(i.e.~away from the constructed component isomorphic to $H$ and these appearances of $K_{3}$).
Thus, the amount of double-counting is at most
$\left( 4|H|-e(H)+1 \right)
\left( \frac{n}{3} \right)^{|H|}
\left(^{m-e(H)-4|H|}_{\phantom{w}3|H|-e(H)}\right)
= \Theta \left( n^{4|H|-e(H)} \right)$,
recalling that $m = \Theta (n)$.

Hence,
we find that the number of \textit{distinct} graphs that we have constructed is at least
$\frac{|\mathcal{G}_{n}| \Theta \left( n^{4|H|-e(H)} \right)}{\Theta \left( n^{4|H|-e(H)} \right)}
= |\mathcal{G}_{n}| \Theta (1)$.
Thus,
recalling that
$|\mathcal{G}_{n}| \geq \frac{\epsilon}{2} |\mathcal{P}(n,m)|$,
we are done.
\phantom{www}
\begin{picture}(1,1)
\put(0,0){\line(1,0){1}}
\put(0,0){\line(0,1){1}}
\put(1,1){\line(-1,0){1}}
\put(1,1){\line(0,-1){1}}
\end{picture} \\
\\

Note that in the previous proof,
we could have constructed a component isomorphic to $H$ directly from an appearance of $H$.
We chose to instead build the component from isolated vertices cut from appearances of $K_{4}$,
as this technique generalises more easily to our next proof,
as we shall now explain.

Recall that when we cut the isolated vertices from the appearances of $K_{4}$,
this involved deleting three edges for each isolated vertex that we created,
which crucially meant that we had enough edges to play with when we wanted to turn these isolated vertices
into a component isomorphic to $H$.
Notice, though, that the proof was only made possible by the fact that we had lots of appearances of $K_{4}$
to choose from,
which was why we needed to restrict $\frac{m}{n}$ to the region $[b,B]$, where $b>1$ and $B<3$.

However, if $e(H) \leq |H|$ then we would have enough edges to play with
even if we only deleted one edge for each isolated vertex that we created.
Thus, we may replace the role of the appearances of $K_{4}$ by pendant edges,
which we know are plentiful even for small values of $\frac{m}{n}$,
by Lemma~\ref{pen5}.
Hence, we may obtain:

\begin{Theorem} \label{gen4}
Let $H$
be a (fixed) connected planar graph with $e(H) \leq |H|$,
let $c>0$ and $B<3$ be constants, and let $m(n) \in [cn,Bn]$ for all $n$.
Then there exist constants $\epsilon (H,c,B) > 0$
and $N(H,c,B)$ such that
\begin{eqnarray*}
& & \mathbf{P}[\textrm{$P_{n,m}$
will have a component isomorphic to $H$} ]
\geq \epsilon \textrm{ for all } n \geq N.
\end{eqnarray*}
\end{Theorem}
\textbf{Proof}
Suppose the result is false.
Then, similarly to with the proof of Theorem~\ref{gen3},
we have a set $\mathcal{G}_{n}$ of at least $\frac{\epsilon}{2} |\mathcal{P}(n,m)|$ graphs with
(i) no components isomorphic to $H$
and (ii) at least $\alpha n$ pendant edges (using Lemma~\ref{pen5}).

\begin{figure} [ht]
\setlength{\unitlength}{0.785cm}
\begin{picture}(20,2.6)(-0.25,0.4)

\put(0.45,0.75){\circle*{0.1}}
\put(1.65,0.75){\circle*{0.1}}
\put(2.85,0.75){\circle*{0.1}}
\put(4.05,0.75){\circle*{0.1}}

\put(2.25,2.25){\oval(5,1.5)}

\put(0.45,0.75){\line(0,1){0.75}}
\put(0.45,1.5){\circle*{0.1}}
\put(1.65,0.75){\line(0,1){0.75}}
\put(1.65,1.5){\circle*{0.1}}
\put(2.85,0.75){\line(0,1){0.75}}
\put(2.85,1.5){\circle*{0.1}}
\put(4.05,0.75){\line(0,1){0.75}}
\put(4.05,1.5){\circle*{0.1}}

\put(6,1.5){\vector(1,0){1}}

\put(10.75,2.25){\oval(5,1.5)}

\put(14,0.75){\line(1,0){0.9}}
\put(14,0.75){\line(3,2){0.45}}
\put(14.45,1.05){\line(0,1){0.45}}
\put(14.9,0.75){\line(-3,2){0.45}}
\put(14,0.75){\circle*{0.1}}
\put(14.9,0.75){\circle*{0.1}}
\put(14.45,1.05){\circle*{0.1}}
\put(14.45,1.5){\circle*{0.1}}

\put(13.8,0.45){\small{$v_{3}$}}
\put(14.8,0.45){\small{$v_{4}$}}
\put(14.55,1.05){\small{$v_{2}$}}
\put(14.55,1.5){\small{$v_{1}$}}

\put(8.95,1.5){\circle*{0.1}}
\put(10.15,1.5){\circle*{0.1}}
\put(11.35,1.5){\circle*{0.1}}
\put(12.55,1.5){\circle*{0.1}}

\put(0.25,0.4){$v_{1}$}
\put(1.45,0.4){$v_{2}$}
\put(2.65,0.4){$v_{3}$}
\put(3.85,0.4){$v_{4}$}

\put(14.25,0){$H$}

\end{picture}
\caption{Constructing a component isomorphic to $H$.} \label{cpt2}

\end{figure}
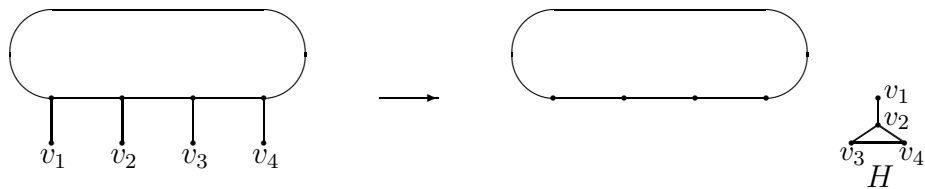

Given a graph $G \in \mathcal{G}_{n}$,
we may delete $|H|$ of the pendant edges
and use the resulting isolated vertices to construct a component isomorphic to $H$ (see Figure~\ref{cpt2}).
If $H$ is a tree, then we should also add one edge in a suitable place somewhere in the rest of the graph.

By similar counting arguments to those used in the proof of Theorem~\ref{gen3},
we achieve our result.
\phantom{qwerty}
\begin{picture}(1,1)
\put(0,0){\line(1,0){1}}
\put(0,0){\line(0,1){1}}
\put(1,1){\line(-1,0){1}}
\put(1,1){\line(0,-1){1}}
\end{picture} \\
\\

By exactly the same proof as for Theorem~\ref{gen4},
using the additional ingredient that add$(n,m) = \omega(n)$ if $\frac{m}{n} \leq 1 + o(1)$
(from Lemma~\ref{add41}),
we also obtain our third result of this section:

\begin{Theorem} \label{tree6}
Let $H$ be a (fixed) tree, let $c>0$ be a constant, and let $m(n) \in [cn, (1+o(1))n]$ as $n \to \infty$.
Then
\begin{displaymath}
\mathbf{P}[P_{n,m} \textrm{ will have a component isomorphic to } H]
\to 1 \textrm{ as } n \to \infty.
\end{displaymath}
\end{Theorem}
\textbf{Proof}
As before,
we suppose that
$\mathbf{P}[P_{n,m}$
will have a component isomorphic to
$H]<1-\epsilon$
for an arbitrary $\epsilon \in (0,1)$,
and that we hence have a set $\mathcal{G}_{n}$
of at least $\frac{\epsilon}{2}|\mathcal{P}(n,m)|$ graphs with
(i) no components isomorphic to $H$
and (ii) at least $\alpha n$ pendant edges.

We proceed as in the proof of Theorem~\ref{gen4},
and the counting is the same except that we now have $\omega(n)$ choices
(instead of just $\Omega(n)$)
for where to add the `extra' edge
after we have constructed the component isomorphic to $H$.
Hence,
we find that we can build $|\mathcal{G}_{n}| \omega(1) = |\mathcal{P}(n,m)| \omega(1)$
distinct graphs in $\mathcal{P}(n,m)$,
which is a contradiction.
\phantom{qwerty}
\begin{picture}(1,1)
\put(0,0){\line(1,0){1}}
\put(0,0){\line(0,1){1}}
\put(1,1){\line(-1,0){1}}
\put(1,1){\line(0,-1){1}}
\end{picture} \\
\\

By using more precise estimates of add$(n,m)$,
lower bounds for the \textit{number} of components in $P_{n,m}$ isomorphic to $H$
can be obtained (see Theorem 38 of~\cite{dow}).
One such result is that $P_{n,m}$ will a.a.s.~have linearly many components isomorphic to any given tree
if we strengthen the upper bound on $\frac{m}{n}$ to $\limsup \frac{m}{n} < 1$, rather than $\limsup \frac{m}{n} \leq 1$.
Since this particular result shall be needed in Section~\ref{sub},
we will now provide a full proof.
The method is exactly the same as with the last two results,
but the equations involved are more complicated:

\begin{Theorem} \label{tree4}
Let $H$ be a (fixed) tree, let $c>0$ and $A<1$ be constants,
and let $m(n) \in [cn, An]$ for all $n$.
Then there exists a constant $\lambda (H,c,A) > 0$ such that
\begin{eqnarray*}
\mathbf{P}[
P_{n,m}
\textrm{ will have less than $\lambda n$
components isomorphic to $H$}]
< e^{- \lambda n } \\
\textrm{for all large $n$}.
\end{eqnarray*}
\end{Theorem}
\textbf{Proof}
By Lemma~\ref{pen5}, we know there exist constants $\alpha >0$, $\beta >0$ and $n_{0}$ such that
$\mathbf{P}[P_{n,m}$ will have less than $\alpha n$ pendant edges] $< e^{- \beta n}$ for all $n \geq n_{0}$.
Let $\lambda$ be a small positive constant
and suppose that there exists a value $n \geq n_{0}$ such that
$\mathbf{P}[P_{n,m}$ will have less than $\lceil \lambda n \rceil$
components isomorphic to $H] \geq e^{-\lceil \lambda n \rceil}$.
Then there is a set $\mathcal{G}_{n}$ of at least a proportion
$e^{- \lambda n} - e^{- \beta n}$ of the graphs in
$\mathcal{P}(n,m)$
with (i) less than $\lambda n$ components isomorphic to $H$ and (ii) at least $\alpha n$ pendant edges.

Without loss of generality,
we may assume that $\lambda$ is small enough
and $n$ large enough
that various inequalities hold during this proof.
In particular,
it is worth noting now that we may assume that
$\alpha n \geq \lceil \lambda n \rceil |H|$ and
$e^{- \lambda n} - e^{- \beta n}
\geq \frac{1}{2}e^{- \lambda n}$.

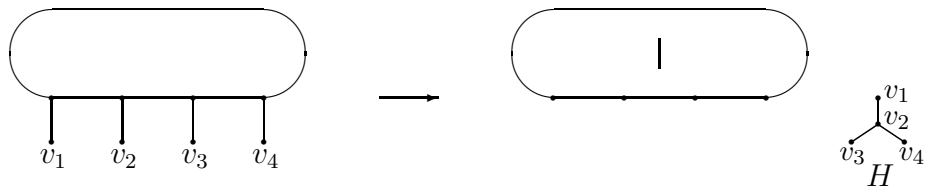
\begin{figure} [ht]
\setlength{\unitlength}{0.785cm}
\begin{picture}(20,2.6)(-0.25,0.4)

\put(0.45,0.75){\circle*{0.1}}
\put(1.65,0.75){\circle*{0.1}}
\put(2.85,0.75){\circle*{0.1}}
\put(4.05,0.75){\circle*{0.1}}

\put(2.25,2.25){\oval(5,1.5)}

\put(0.45,0.75){\line(0,1){0.75}}
\put(0.45,1.5){\circle*{0.1}}
\put(1.65,0.75){\line(0,1){0.75}}
\put(1.65,1.5){\circle*{0.1}}
\put(2.85,0.75){\line(0,1){0.75}}
\put(2.85,1.5){\circle*{0.1}}
\put(4.05,0.75){\line(0,1){0.75}}
\put(4.05,1.5){\circle*{0.1}}

\put(6,1.5){\vector(1,0){1}}

\put(10.75,2.25){\oval(5,1.5)}

\put(14,0.75){\line(3,2){0.45}}
\put(14.45,1.05){\line(0,1){0.45}}
\put(14.9,0.75){\line(-3,2){0.45}}
\put(14,0.75){\circle*{0.1}}
\put(14.9,0.75){\circle*{0.1}}
\put(14.45,1.05){\circle*{0.1}}
\put(14.45,1.5){\circle*{0.1}}

\put(13.8,0.45){\small{$v_{3}$}}
\put(14.8,0.45){\small{$v_{4}$}}
\put(14.55,1.05){\small{$v_{2}$}}
\put(14.55,1.5){\small{$v_{1}$}}

\put(8.95,1.5){\circle*{0.1}}
\put(10.15,1.5){\circle*{0.1}}
\put(11.35,1.5){\circle*{0.1}}
\put(12.55,1.5){\circle*{0.1}}

\put(0.25,0.4){$v_{1}$}
\put(1.45,0.4){$v_{2}$}
\put(2.65,0.4){$v_{3}$}
\put(3.85,0.4){$v_{4}$}

\put(10.75,2){\line(0,1){0.5}}

\put(14.25,0){$H$}

\end{picture}
\caption{Constructing a component isomorphic to $H$.} \label{cpt3}

\end{figure}

To build graphs with at least $\lambda n$ components isomorphic to $H$,
one can start with a graph $G \in \mathcal{G}_{n}$ ($|\mathcal{G}_{n}|$ choices),
delete $\lceil \lambda n \rceil |H|$ of the pendant edges
$\left( \textrm{ at least } \left(^{\phantom{w}\lceil \alpha n \rceil}
_{\lceil \lambda n \rceil |H|} \right) \textrm{ choices } \right)$,
and insert edges between $\lceil \lambda n \rceil |H|$ of the newly-isolated vertices
(choosing one from each pendant edge)
in such a way that they now form $\lceil \lambda n \rceil$ components isomorphic to $H$
$\left( \textrm{at least }
\left( ^{\phantom{i}\lceil \lambda n \rceil |H|} _{|H|, \ldots, |H|} \right)
\frac{1}{\lceil \lambda n \rceil !} \textrm{ choices} \right)$.
We should then add $\lceil \lambda n \rceil$ edges somewhere in the rest of the graph
(i.e.~away from our newly constructed components)
to maintain the correct number of edges overall
$\left( \textrm{we have at least }
\prod_{i=0}^{\lceil \lambda n \rceil-1}
\textrm{add}(n-\lceil \lambda n \rceil |H|,m-\lceil \lambda n \rceil |H|+i) \geq \right.$
$(\textrm{add}(n-\lceil \lambda n \rceil |H|,m))^{\lceil \lambda n \rceil}
\textrm{ choices for this} \Big)$.
See Figure~\ref{cpt3}.

Hence, the number of ways that we have to build
(not necessarily distinct) graphs in $\mathcal{P}(n,m)$ that have
at least $\lambda n$ components isomorphic to $H$
is at least
\begin{eqnarray*}
& & \frac{ \lceil \alpha n \rceil! }
{ \left( \lceil \lambda n \rceil |H| \right)! \left( \lceil \alpha n \rceil - \lceil \lambda n \rceil |H| \right)! }
\frac{ \left( \lceil \lambda n \rceil |H| \right)! }
{ \left( |H|! \right)^{\lceil \lambda n \rceil} }
\frac{1}{\lceil \lambda n \rceil!} \\
& & \cdot
(\textrm{add}(n-\lceil \lambda n \rceil |H|,m))^{\lceil \lambda n \rceil} |\mathcal{G}_{n}| \\
& \geq & \left( \lceil \alpha n \rceil - \lceil \lambda n \rceil |H| \right)^{\lceil \lambda n \rceil |H|}
\left( \frac{1}{|H|!} \right)^{\lceil \lambda n \rceil}
\frac{1}{\lceil \lambda n \rceil!} \\
& & \cdot
(\textrm{add}(n-\lceil \lambda n \rceil |H|,m))^{\lceil \lambda n \rceil} |\mathcal{G}_{n}| \\
& \geq & \left( \left( \frac{\alpha n}{2} \right)^{|H|}
\left( \frac{1}{|H|!} \right)
(\textrm{add}(n-\lceil \lambda n \rceil |H|,m)) \right)^{\lceil \lambda n \rceil}
\frac{1}{\lceil \lambda n \rceil!} |\mathcal{G}_{n}| \\
& & \textrm{(since we may assume that $\lambda$ is sufficiently small
and $n$ sufficiently large} \\
& & \textrm{that $\lceil \alpha n \rceil - \lceil \lambda n \rceil |H| \geq \frac{\alpha n}{2}$)}.
\end{eqnarray*}

Let us now consider the amount of double-counting: \\
Each of our constructed graphs will contain at most $\lceil \lambda n \rceil (|H|+3)-1$ components isomorphic to $H$
(since there were at most $\lceil \lambda n \rceil -1$ already in $G$;
we have deliberately added $\lceil \lambda n \rceil$;
and we may have created at most one extra one each time
we deleted a pendant edge or added an edge in the rest of the graph),
so we have at most
$\left(^{\lceil \lambda n \rceil (|H|+3)-1}_{\phantom{www}\lceil \lambda n \rceil}\right) \leq
\frac{1}{\lceil \lambda n \rceil !} (\lceil \lambda n \rceil (|H|+3))^{\lceil \lambda n \rceil}$
possibilities for which are our created components.
We then have at most $n^{\lceil \lambda n \rceil |H|}$
possibilities for where the vertices in our created components were attached originally
and at most $\left(^{\phantom{q}m}_{\lceil \lambda n \rceil}\right) \leq (3n)^{\lceil \lambda n \rceil}$
possibilities for which edges was added.
Thus, the amount of double-counting is at most
$\frac{1}{\lceil \lambda n \rceil!} \left( \lceil \lambda n \rceil (|H|+3)n^{|H|}3n \right)^{\lceil \lambda n \rceil}$.

Hence, putting everything together,
we find that the number of \textit{distinct} graphs in $\mathcal{P}(n,m)$ that have at least
$\lambda n$ components
isomorphic to $H$ is at least
\begin{displaymath}
\left( \left( \frac{\alpha}{2} \right)^{|H|} \frac{1}{|H|!} (\textrm{add}(n- \lceil \lambda n \rceil |H|,m))
\frac{1}{ \lceil \lambda n \rceil (|H|+3)3n } \right)^{\lceil \lambda n \rceil} |\mathcal{G}_{n}|.
\end{displaymath}

Recall that $m \leq An$,
where $A<1$.
Thus,
we may assume that $\lambda$ is sufficiently small
and $n$ sufficiently large that
$m \leq \frac{A+1}{2}(n-\lceil \lambda n \rceil |H|)$.
Hence,
by Lemma~\ref{add2},
we have
add$(n- \lceil \lambda n \rceil |H|,m) \geq
(1+o(1)) \left( \frac{ \left( 1- \frac{A+1}{2} \right)^{2} }{2} \right) n^{2}
= (1+o(1)) \frac{(1-A)^{2}}{8} n^{2}$.

Therefore, we find that
the number of graphs in $\mathcal{P}(n,m)$ that have at least $\lambda n$ components
isomorphic to $H$ is at least
\begin{displaymath}
\left( (1+o(1)) \frac{\alpha^{|H|}}{2^{|H|}|H|!3(|H|+3)}
\left( \frac{ (1-A)^{2} }{ 8 \lambda } \right)
\right)^{\lceil \lambda n \rceil} |\mathcal{G}_{n}|.
\end{displaymath}
But this is more than $|\mathcal{P}(n,m)|$ for large $n$,
if $\lambda$ is sufficiently small,
since we recall that
$|\mathcal{G}_{n}| \geq \frac{1}{2} e^{- \lambda n} |\mathcal{P}(n,m)|$.
Thus, by proof by contradiction, it must be that
$\mathbf{P}[P_{n,m}$ will have less than $\lambda n$
components isomorphic to $H] < e^{- \lambda n}$ for all large $n$.
\phantom{qwerty}
\setlength{\unitlength}{.25cm}
\begin{picture}(1,1)
\put(0,0){\line(1,0){1}}
\put(0,0){\line(0,1){1}}
\put(1,1){\line(-1,0){1}}
\put(1,1){\line(0,-1){1}}
\end{picture} \\

\section{Components II: Upper Bounds} \label{cyccpt}

In this section,
we shall produce upper bounds for
$\mathbf{P} := \mathbf{P}[P_{n,m}$ will have a component isomorphic to $H]$
to complement the lower bounds of Section~\ref{cptlow}.

We will start with the case $0 < \liminf \frac{m}{n} \leq \limsup \frac{m}{n} \leq 1$,
for which we have seen $\mathbf{P} \to 1$ if $H$ is a tree
and $\liminf \mathbf{P} > 0$ if $H$ is unicyclic.
In this section,
we shall complete matters by showing $\mathbf{P} \to 0$
if $H$ is multicyclic (see Theorem~\ref{cyc41})
and $\limsup \mathbf{P} < 1$ if $H$ is unicyclic (see Theorem~\ref{cyc52}).

We will then deal with the case when $\liminf \frac{m}{n} > 1$,
for which we have seen $\liminf \mathbf{P} > 0$
for all connected planar $H$
if we also have $\limsup \frac{m}{n} < 3$.
By examining the probability that $P_{n,m}$ is connected,
we will now show $\mathbf{P} \to 0$ if $\frac{m}{n} \to 3$ (see Theorem~\ref{conn5})
and $\limsup \mathbf{P} < 1$ if $\liminf \frac{m}{n} > 1$ (see Theorem~\ref{conn4}). \\

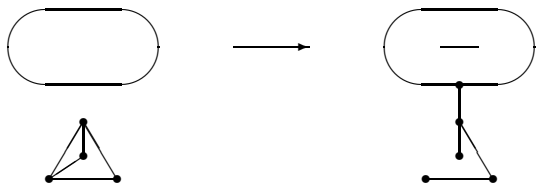
\begin{figure} [ht]
\setlength{\unitlength}{1cm}
\begin{picture}(20,2.25)

\put(3.05,0){\line(1,0){0.9}}
\put(3.05,0){\line(3,2){0.45}}
\put(3.05,0){\line(3,5){0.45}}
\put(3.5,0.3){\line(0,1){0.45}}
\put(3.95,0){\line(-3,5){0.45}}
\put(3.05,0){\circle*{0.1}}
\put(3.95,0){\circle*{0.1}}
\put(3.5,0.3){\circle*{0.1}}
\put(3.5,0.75){\circle*{0.1}}

\put(3.5,1.75){\oval(2,1)}

\put(5.5,1.75){\vector(1,0){1}}

\put(8.05,0){\line(1,0){0.9}}
\put(8.5,0.3){\line(0,1){0.45}}
\put(8.95,0){\line(-3,5){0.45}}
\put(8.05,0){\circle*{0.1}}
\put(8.95,0){\circle*{0.1}}
\put(8.5,0.3){\circle*{0.1}}
\put(8.5,0.75){\circle*{0.1}}

\put(8.5,1.75){\oval(2,1)}

\put(8.25,1.75){\line(1,0){0.5}}

\put(8.5,0.75){\line(0,1){0.5}}
\put(8.5,1.25){\circle*{0.1}}

\end{picture}

\caption{Redistributing edges from our multicyclic component.} \label{cpt4}
\end{figure}

We start with our aforementioned result for multicyclic components when $\frac{m}{n} \leq 1+o(1)$:

\begin{Theorem} \label{cyc41}
Let $H$ be a (fixed) multicyclic connected planar graph
and let $m(n) \leq (1+o(1))n$.
Then
\begin{displaymath}
\mathbf{P}[P_{n,m} \textrm{ will have a component isomorphic to } H] \to 0 \textrm{ as } n \to \infty.
\end{displaymath}
\end{Theorem}
\textbf{Proof}
Let $\mathcal{G}_{n}$ denote the set of graphs in $\mathcal{P}(n,m)$ with a component isomorphic to $H$.
For each graph $G \in \mathcal{G}_{n}$,
let us delete $2$ edges from a component $H^{\prime} (=H^{\prime}_{G})$ isomorphic to $H$
in such a way that we do not disconnect the component.
Let us then insert one edge between a vertex in the remaining component and a vertex elsewhere in the graph.
We have $|H|(n-|H|)$ ways to do this, and planarity is maintained.
Let us then also insert one other edge into the graph,
without violating planarity
(see Figure~\ref{cpt4}).
We have at least $(\textrm{add}(n,m))=\omega(n)$
choices for where to place this second edge,
by Lemma~\ref{add41}.
Thus, we can construct $|\mathcal{G}_{n}| \omega \left( n^{2} \right)$
(not necessarily distinct) graphs in $\mathcal{P}(n,m)$.

Given one of our constructed graphs, there are $m=O(n)$
possibilities for the edge that was inserted last.
There are then at most $m-1=O(n)$ possibilities for the other edge that was inserted.
Since one of the two vertices incident with this edge must belong to $V(H^{\prime})$,
we then have at most two possibilities for $V(H^{\prime})$
and then at most $\left( ^{\left( ^{|H|}_{\phantom{q}2} \right)}_{\phantom{qq}2} \right) = O(1)$
possibilities for $E(H^{\prime})$.
Thus, we have built each graph at most $O \left( n^{2} \right)$ times,
and so
$\frac{|\mathcal{G}_{n}|}{|\mathcal{P}(n,m)|} = \frac{O \left( n^{2} \right)}{\omega \left( n^{2} \right)} \to 0$.
\phantom{qqqqi}
\setlength{\unitlength}{.25cm}
\begin{picture}(1,1)
\put(0,0){\line(1,0){1}}
\put(0,0){\line(0,1){1}}
\put(1,1){\line(-1,0){1}}
\put(1,1){\line(0,-1){1}}
\end{picture} \\
\\

We shall now look at unicyclic components.
The basic argument will be the same as that of Theorem~\ref{cyc41},
i.e.~we will start with $\mathcal{G}_{n}$,
the set of graphs in $\mathcal{P}(n,m)$ with a component isomorphic to $H$,
and redistribute edges from such components to construct other graphs in $\mathcal{P}(n,m)$.
This time we will only be able to delete one edge from each component,
and so we will only be able to show
$\frac{|\mathcal{G}_{n}|}{|\mathcal{P}(n,m)|} = O(1)$,
rather than $o(1)$.
Hence,
it will be crucial to keep track of the constants involved in the calculations,
so that we can try to show $\limsup_{n \to \infty} \frac{|\mathcal{G}_{n}|}{|\mathcal{P}(n,m)|} < 1$.
Unfortunately,
it turns out that the constants are actually only small enough for certain $H$,
such as when $H$ is a cycle.
However,
it is fairly simple to relate the probability that a given component is isomorphic to $H$
to the probability that it is a cycle,
and so we may deduce that the result actually holds for any unicyclic $H$.

\begin{figure} [ht]
\setlength{\unitlength}{1cm}
\begin{picture}(20,2.25)

\put(3.05,0){\line(1,0){0.9}}
\put(3.05,0){\line(3,5){0.45}}
\put(3.95,0){\line(-3,5){0.45}}
\put(3.05,0){\circle*{0.1}}
\put(3.95,0){\circle*{0.1}}
\put(3.5,0.75){\circle*{0.1}}

\put(3.5,1.75){\oval(2,1)}

\put(5.5,1.75){\vector(1,0){1}}

\put(8.05,0){\line(1,0){0.9}}
\put(8.5,0.75){\line(0,1){0.5}}
\put(8.95,0){\line(-3,5){0.45}}
\put(8.05,0){\circle*{0.1}}
\put(8.95,0){\circle*{0.1}}
\put(8.5,0.75){\circle*{0.1}}

\put(8.5,1.75){\oval(2,1)}

\put(8.5,1.25){\circle*{0.1}}

\end{picture}

\caption{Redistributing an edge from our cycle.} \label{subcpt4}
\end{figure}
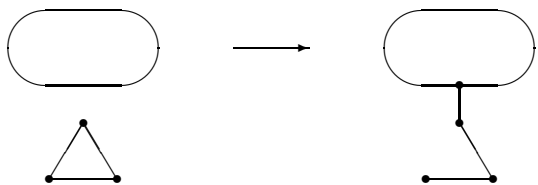

\begin{Theorem} \label{cyc52}
Let $H$ be a (fixed) unicyclic connected planar graph.
Then, given any $m(n)$,
\begin{displaymath}
\limsup_{n \to \infty}
\mathbf{P}[P_{n,m} \textrm{ will have a component isomorphic to } H] < 1.
\end{displaymath}
\end{Theorem}
\textbf{Proof}
Let $\mathcal{G}_{n}$ denote the set of graphs in $\mathcal{P}(n,m)$ with a component
that is a \textit{cycle} of order $|H|$.
For each graph $G \in \mathcal{G}_{n}$,
let us delete an edge from such a component,
leaving a spanning path.
Note that we have $|H|$ choices for this edge.
Let us then insert an edge between a vertex in the spanning path and a vertex elsewhere
(see Figure~\ref{subcpt4}).
We have $|H|(n-|H|)$ ways to do this and planarity is maintained.
Thus, we have constructed $|\mathcal{G}_{n}||H|^{2}(n-|H|)$ graphs in $\mathcal{P}(n,m)$.

Given one of our constructed graphs,
there are at most $n-1$ possibilities for which is the inserted edge,
since it must be a cut-edge.
There are then at most $2$ possibilities for the vertex set of our original cycle
(at most one possibility for each endpoint of the edge),
and then only one possibility for where the inserted edge was originally.
Thus, we have built each graph at most $2(n-1)$ times.
Therefore, $|\mathcal{G}_{n}| \leq \frac{2(n-1)}{|H|^{2}(n-|H|)} |\mathcal{P}(n,m)|$,
and hence there exists an $\epsilon > 0$
such that $|\mathcal{G}_{n}^{c}| > \epsilon |\mathcal{P}(n,m)|$
for all $n$,
where $\mathcal{G}_{n}^{c}$ denotes $\mathcal{P}(n,m) \setminus \mathcal{G}_{n}$.

If $H$ is itself a cycle,
then we are done.
If not,
let $r$ satisfy $\frac{|H|!}{r} \leq \frac{\epsilon}{2}$
and let $\mathcal{J}_{n}$ denote the set of graphs in $\mathcal{P}(n,m)$ with
more than $r$ components isomorphic to $H$.
Let us now consider the set $\mathcal{J}_{n} \cap \mathcal{G}_{n}^{c}$,
i.e.~the set of graphs in $\mathcal{P}(n,m)$ with
more than $r$ components isomorphic to $H$,
but no components that are cycles of order $|H|$.

For each graph $J \in \mathcal{J}_{n} \cap \mathcal{G}_{n}^{c}$,
delete a component isomorphic to $H$ (we have more than $r$ choices for this)
and replace it with a cycle of order $|H|$.
Note that each graph is built at most $|H|!$ times,
since we will be able to see exactly where the modified component is,
and so we can construct at least
$\frac{r}{|H|!} |\mathcal{J}_{n} \cap \mathcal{G}_{n}^{c}|$
distinct graphs in $\mathcal{P}(n,m)$.
Hence,
$|\mathcal{J}_{n} \cap \mathcal{G}_{n}^{c}| \leq \frac{|H|!}{r} |\mathcal{P}(n,m)|
\leq \frac{\epsilon}{2} |\mathcal{P}(n,m)|$.

Finally,
let $\mathcal{J}_{n}^{c}$ denote $\mathcal{P}(n,m) \setminus \mathcal{J}_{n}$
and let us consider the set $\mathcal{J}_{n}^{c} \cap \mathcal{G}_{n}^{c}$
(i.e.~the set of graphs in $\mathcal{P}(n,m)$ with
at most $r$ components isomorphic to $H$
and no components that are cycles of order $|H|$).
Note that we have
$|\mathcal{J}_{n}^{c} \cap \mathcal{G}_{n}^{c}|
=|\mathcal{G}_{n}^{c}|-|\mathcal{J}_{n} \cap \mathcal{G}_{n}^{c}|
\geq \frac{\epsilon}{2}|\mathcal{P}(n,m)|$,
by the previous paragraph.

For each graph $L \in \mathcal{J}_{n}^{c} \cap \mathcal{G}_{n}^{c}$,
let us change each component isomorphic to $H$ into a cycle.
The amount of double-counting will be at most $\left( |H|! \right)^{r}$,
since we will be able to see exactly where the modified components are,
and so we can construct at least
$\frac{|\mathcal{J}_{n}^{c} \cap \mathcal{G}_{n}^{c}|}{\left( |H|! \right)^{r}}
\geq \frac{\epsilon |\mathcal{P}(n,m)|}{2\left( |H|! \right)^{r}}$
distinct graphs in $\mathcal{P}(n,m)$ without any components isomorphic to $H$.
\phantom{qwerty}
\setlength{\unitlength}{.25cm}
\begin{picture}(1,1)
\put(0,0){\line(1,0){1}}
\put(0,0){\line(0,1){1}}
\put(1,1){\line(-1,0){1}}
\put(1,1){\line(0,-1){1}}
\end{picture} \\
\\

As mentioned at the start of this section,
we shall now obtain upper bounds for
$\mathbf{P}[P_{n,m}$ will have a component isomorphic to $H]$
when $\liminf \frac{m}{n} > 1$
by investigating the probability that $P_{n,m}$ will be connected,
in which case it clearly won't contain any components of order less than $n$ at all.

We start with the case when $\frac{m}{n} \to 3$:

\begin{Theorem} \label{conn5}
Let $m(n)=3n-o(n)$.
Then
\begin{displaymath}
\mathbf{P}[P_{n,m} \textrm{ will be connected}] \to 1 \textrm{ as } n \to \infty.
\end{displaymath}
\end{Theorem}
\textbf{Proof}
Let $G \in \mathcal{P}(n,m)$ and let us consider how many triangles in $G$ contain at least one vertex with degree $\leq 6$.
We shall call such triangles `good' triangles.

First, note that (assuming $n \geq 3$)
$G$ may be extended to a triangulation by inserting $3n-6-m = o(n)$ `phantom' edges.
Let $d_{i}$ denote the number of vertices of degree $i$ in such a triangulation.
Then $7 \sum_{i \geq 7} d_{i} \leq \sum_{i \geq 1} i d_{i} = 2(3n-6)$.
Thus, $\sum_{i \geq 7} d_{i} < \frac{6n}{7}$ and so $\sum_{i \leq 6} d_{i} > \frac{n}{7}$.

For $n>3$, each of these vertices of small degree will be in at least three faces of the triangulation,
all of which will be good triangles.
This counts each good triangle at most three times
(once for each vertex),
so our triangulation must have at least
$3 \cdot \frac{n}{7} \cdot \frac{1}{3} = \frac{n}{7}$ good triangles that are \textit{faces}.
Each of our phantom edges is in exactly two faces of the triangulation,
so our original graph $G$ must contain at least $\frac{n}{7} + o(n)$ of our good triangles
(note that these triangles will still be `good',
since the degrees of the vertices will be at most what they were in the triangulation).

We will now consider how many cut-edges a graph in $\mathcal{P}(n,m)$ may have.
If we delete all $c$ cut-edges,
then the remaining graph will consist of $b$, say, blocks,
each of which is either $2$-edge-connected or is an isolated vertex.
Note that the graph formed by condensing each block to a single node and re-inserting the cut-edges must be acyclic,
so $c \leq b-1$.
Label the blocks $1,2, \ldots b$ and let $n_{i}$ denote the number of vertices in block $i$.
Then the number of edges in block $i$ is at most $3n_{i}-6$ if $n_{i} \geq 3$ and is $0=3n_{i}-3$ otherwise
(since $n_{i} < 3$ implies that $n_{i} = 1$).
Thus, $m \leq \sum_{i=1}^{b} (3n_{i}-3) + c = 3n-3b+c < 3n-2c$,
and so $c < \frac{3n-m}{2} = o(n)$.

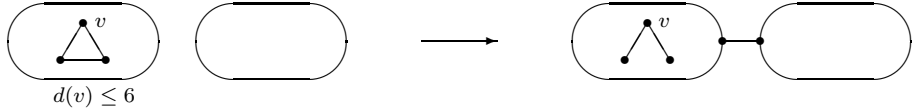
\begin{figure} [ht]
\setlength{\unitlength}{1cm}
\begin{picture}(20,1.3)(0,-0.3)

\put(0.7,0.25){\line(1,0){0.6}}
\put(0.7,0.25){\line(3,5){0.3}}
\put(1.3,0.25){\line(-3,5){0.3}}
\put(0.7,0.25){\circle*{0.1}}
\put(1.3,0.25){\circle*{0.1}}
\put(1,0.75){\circle*{0.1}}
\put(1.15,0.7){\scriptsize{$v$}}
\put(0.6,-0.3){\scriptsize{$d(v) \leq 6$}}

\put(1,0.5){\oval(2,1)}
\put(3.5,0.5){\oval(2,1)}

\put(5.5,0.5){\vector(1,0){1}}

\put(8.2,0.25){\line(3,5){0.3}}
\put(8.8,0.25){\line(-3,5){0.3}}
\put(8.2,0.25){\circle*{0.1}}
\put(8.8,0.25){\circle*{0.1}}
\put(8.5,0.75){\circle*{0.1}}
\put(8.65,0.7){\scriptsize{$v$}}

\put(8.5,0.5){\oval(2,1)}
\put(11,0.5){\oval(2,1)}

\put(9.5,0.5){\line(1,0){0.5}}
\put(9.5,0.5){\circle*{0.1}}
\put(10,0.5){\circle*{0.1}}

\end{picture}

\caption{Constructing our new graph.} \label{cpt5}
\end{figure}

We now come to the main part of the proof.
Let $\mathcal{G}_{n}$ denote the set of graphs in $\mathcal{P}(n,m)$ that are not connected,
and choose a graph $G \in \mathcal{G}_{n}$.
Choose a good triangle in $G$ (at least $\frac{n}{7} + o(n)$ choices)
and delete an edge that is opposite a vertex with degree $\leq 6$.
Then insert an edge between two vertices in different components
(we have $a$, say, choices for this edge).
See Figure~\ref{cpt5}.

As mentioned in the previous proof,
the number of possible edges between disjoint sets $X$ and $Y$ is $|X||Y|$ and
if $|X| \leq |Y|$ then $|X||Y| > (|X|-1)(|Y|+1)$,
so it follows that the number of choices for the edge to insert is minimized when
we have one isolated vertex and
one component of $n-1$ vertices.
Thus, $a \geq n-1$ and so we have created at least
$|\mathcal{G}_{n}| \left( \frac{n}{7} - o(n) \right) (n-1)$
graphs in $\mathcal{P}(n,m)$.

Given one of our created graphs,
there are at most $o(n)$ possibilities for which edge was inserted,
since it must be a cut-edge.
There are then at most $\left( ^{6}_{2} \right)n$ possibilities for where the deleted edge was originally,
since it must have been between two neighbours of a vertex with degree $\leq 6$
(we have at most $n$ possibilities for this vertex and then
at most $\left( ^{6}_{2} \right)$ possibilities for its neighbours).
Hence, we have built each graph at most $o(n^{2})$ times.

Thus,
$|\mathcal{P}(n,m)| \geq
\frac{(\frac{1}{7})n^{2} + o \left( n^{2} \right)}{ o \left( n^{2} \right)}
|\mathcal{G}_{n}|$,
and so
$\frac{|\mathcal{G}_{n}|}{|\mathcal{P}(n,m)|} \leq \frac{o \left( n^{2} \right)}{\Theta \left( n^{2} \right)}
\to 0$.
\phantom{qwerty}
\setlength{\unitlength}{.25cm}
\begin{picture}(1,1)
\put(0,0){\line(1,0){1}}
\put(0,0){\line(0,1){1}}
\put(1,1){\line(-1,0){1}}
\put(1,1){\line(0,-1){1}}
\end{picture} \\
\\

It now only remains to look at when
$1 < \liminf \frac{m}{n} \leq \limsup \frac{m}{n} < 3$.

In~\cite{gim},
generating function techniques were used to produce rather precise asymptotic estimates for both
$|\mathcal{P}(n, \lfloor qn \rfloor)|$ and $|\mathcal{P}_{c}(n, \lfloor qn \rfloor)|$,
the number of connected graphs in $\mathcal{P}(n, \lfloor qn \rfloor)$,
when $q \in (1,3)$:

\begin{Proposition}[\cite{gim}, Theorem 3] \label{gim T3}
Let $q \in (1,3)$ be a constant.
Then
\begin{eqnarray*}
|\mathcal{P}(n,\lfloor qn \rfloor)| & \sim &
g(q) \cdot \left( u(q) \right) ^{qn - \lfloor qn \rfloor} \cdot n^{-4} \gamma (q)^{n}n! \\
\textrm{and } \phantom{w} |\mathcal{P}_{c}(n,\lfloor qn \rfloor)| & \sim &
g_{c}(q) \cdot \left( u(q) \right) ^{qn - \lfloor qn \rfloor} \cdot n^{-4} \gamma (q)^{n}n!,
\end{eqnarray*}
where $g(q)$, $g_{c}(q)$, $u(q)$ and $\gamma (q)$ are computable analytic functions.
\end{Proposition}

Clearly,
Proposition~\ref{gim T3} may be used to obtain the exact asymptotic limit for
$\mathbf{P}[P_{n, \lfloor qn \rfloor}$ will be connected] when $q \in (1,3)$,
and it turns out that this limit is strictly greater than $0$ for all $q > 1$.
It is no surprise that this result holds uniformly,
and a combinatorial proof is given in~\cite{dow}.
Hence,
we have

\begin{Theorem}[\cite{dow}, Theorem 44] \label{conn4}
Let $b>1$ be a constant and let $m(n) \in [bn, 3n-6]$.
Then there exists a constant $c(b) > 0$ such that
\begin{displaymath}
\mathbf{P}[P_{n,m} \textrm{ will be connected}] > c \textrm{ for all } n.
\end{displaymath}
\end{Theorem}

\phantom{p}

Note that we now have a complete description of
$\mathbf{P}[P_{n,m}$ will have a component isomorphic to $H]$,
in terms of exactly when it is/isn't bounded away from $0$ and/or $1$.
By combining Theorems~\ref{gen3} and~\ref{tree6} with Theorem~\ref{conn4},
we also have a complete picture of
$\mathbf{P}[P_{n,m}$ will be connected]. \\

\section{Subgraphs} \label{sub}

We have now finished looking at the probability that
$P_{n,m}$ will have a component isomorphic to $H$.
In this section,
we will instead investigate
$\mathbf{P}^{\prime} := \mathbf{P}[P_{n,m}$ will have a \textit{subgraph} isomorphic to $H]$.
We already know from Lemma~\ref{gen2} that
$\mathbf{P}^{\prime} \to 1$ for all $H$ if $1 < \liminf \frac{m}{n} \leq \limsup \frac{m}{n} < 3$,
and we shall soon see (in Theorem~\ref{sub4}) that we can actually drop the $\limsup \frac{m}{n} < 3$ condition.
Hence, the interest lies with the case when $0 < \liminf \frac{m}{n} \leq \limsup \frac{m}{n} \leq 1$.

It follows from our results on components that $\mathbf{P}^{\prime} \to 1$
for this region if $H$ is a tree
and that $\limsup \mathbf{P}^{\prime} > 0$ if $H$ is a connected unicyclic graph.
In this section,
we shall see for the latter case that $\limsup \mathbf{P}^{\prime} < 1$ if $\limsup \frac{m}{n} < 1$
(see Theorem~\ref{cyc54}),
but that $\mathbf{P}^{\prime} \to 1$ if $\frac{m}{n} \to 1$ (see Theorem~\ref{unisub1}).
If $H$ is a connected multicyclic graph,
we shall see (in Theorem~\ref{msub3}) that $\mathbf{P}^{\prime} \to 0$ if $\limsup \frac{m}{n} < 1$.
This leaves the remaining unanswered question of what happens when $H$ is multicyclic and $\frac{m}{n} \to 1$. \\

We now start with a short discussion of the case when $\frac{m}{n} \to 3$.
Since we know from Lemma~\ref{gen2} that
$\mathbf{P}[P_{n,m}$ will have a copy of $H] \to 1$ as $n \to \infty$ for all planar $H$
if $1 < \liminf \frac{m}{n} \leq \limsup \frac{m}{n} < 3$,
it is intuitive that the result ought to also hold when $\frac{m}{n} \to 3$.
Clearly,
we wouldn't expect to have any appearances of $H$ if $m$ is close to $3n-6$,
since appearances involve cut-edges,
but it turns out that the same proof does work if we replace appearances of $H$
by `6-appearances' of triangulations containing $H$,
where a 6-appearance is similar to an appearance
but with six edges (from three vertices) connecting it to the rest of the graph instead of just one
(see~\cite{dow} for details).

Hence, we obtain:

\begin{Theorem}[\cite{dow}, Theorem 61] \label{sub4}
Let $H$ be a (fixed) planar graph and let $m(n)$ satisfy $\liminf_{n \to \infty} \frac{m}{n} > 1$.
Then there exists a constant $\alpha > 0$ such that
\begin{displaymath}
\mathbf{P}[\textrm{$P_{n,m}$
will \emph{not} have a set of at least $\alpha n$ copies of $H$}]
= e^{- \Omega (n)}.
\end{displaymath}
\end{Theorem}

\phantom{p}

We will now spend the remainder of this paper looking at the case when we have $\frac{m}{n} \leq 1+o(1)$.
We start with an upper bound for unicyclic graphs when $\limsup \frac{m}{n} < 1$:

\begin{Theorem} \label{cyc54}
Let $H$ be a (fixed) unicyclic connected planar graph,
let $A<1$ be a constant,
and let $m(n) \leq An$ for all $n$.
Then
\begin{displaymath}
\limsup_{n \to \infty}
\mathbf{P}[P_{n,m} \textrm{ will have a copy of }H] < 1.
\end{displaymath}
\end{Theorem}
\textbf{Sketch of Proof}
Let $C$ denote the unique cycle in $H$.
Clearly,
it suffices to show
$\limsup_{n \to \infty}
\mathbf{P}[P_{n,m} \textrm{ will have a copy of }C] < 1.$

By Theorems~\ref{tree4} and~\ref{cyc52},
we know that there must be a decent proportion of graphs
with lots of isolated vertices,
but no components isomorphic to $C$.
If such a graph contains lots of copies of $C$,
then we may `transfer' the edges of one of these to some isolated vertices to construct a component isomorphic to $C$,
and the amount of double-counting will be small since
there were no such components originally.
Hence, we can use this idea to show that the proportion of graphs
with lots of copies of $C$, lots of isolated vertices and no components isomorphic to $C$ must be small,
and so there must be a decent proportion of graphs with few copies of $C$
(and lots of isolated vertices and no components isomorphic to $C$).

Given such a graph,
we may then destroy the few copies of $C$
by deleting the edges from them and inserting edges between components of the graph. \\
\\
\textbf{Full Proof}
As mentioned in the sketch-proof,
it clearly suffices for us to show
$\limsup_{n \to \infty}
\mathbf{P}[P_{n,m} \textrm{ will have a copy of }C] < 1$,
where $C$ denotes the unique cycle in $H$.

Let $b \in \left( 0, \frac{1}{2} \right)$ be a fixed constant and suppose $m \in [bn,An]$ for all $n$.
Then, by Theorem~\ref{tree4} (with $H$ as an isolated vertex),
there exists a constant $\beta > 0$ such that
$\mathbf{P}[P_{n,m} \textrm{ will have at least $\beta n$ isolated vertices}] \to 1 \textrm{ as } n \to \infty$.
Clearly, $P_{n,m}$ must also have at least $n-2bn$ isolated vertices if $m \leq bn$,
so by setting $\alpha = \min \{ \beta, 1-2b \}$
we find that
$\mathbf{P}[P_{n,m} \textrm{ will have at least $\alpha n$ isolated vertices}] \to 1$ as $n \to \infty$
whenever $m \leq An$.

Let $\mathcal{G}_{n}$ denote the set of graphs in $\mathcal{P}(n,m)$ with at least $\alpha n$ isolated vertices
and no components isomorphic to $C$.
Then, using Theorem~\ref{cyc52},
we now know that there exists an $\epsilon > 0$ such that
$|\mathcal{G}_{n}| \geq \epsilon |\mathcal{P}(n,m)|$
for all sufficiently large $n$.

Let $r$ satisfy $\frac{(|C|+1)!}{r \alpha^{|C|}} < \frac{\epsilon}{2}$
and let $\mathcal{J}_{n}$ denote the set of graphs in $\mathcal{P}(n,m)$
whose maximal number of edge-disjoint copies of $C$ is more than $r$.
Let us now consider the set $\mathcal{J}_{n} \cap \mathcal{G}_{n}$,
i.e.~the set of graphs in $\mathcal{P}(n,m)$ with at least $\alpha n$ isolated vertices,
more than $r$ edge-disjoint copies of $C$,
but no components isomorphic to $C$.

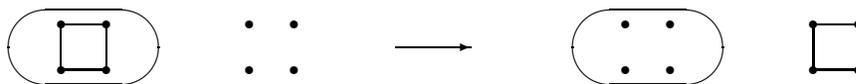
\begin{figure} [ht]
\setlength{\unitlength}{1cm}
\begin{picture}(20,1.3)(-0.35,-0.3)

\put(0.7,0.2){\line(1,0){0.6}}
\put(0.7,0.2){\line(0,1){0.6}}
\put(0.7,0.8){\line(1,0){0.6}}
\put(1.3,0.2){\line(0,1){0.6}}
\put(0.7,0.2){\circle*{0.1}}
\put(1.3,0.2){\circle*{0.1}}
\put(0.7,0.8){\circle*{0.1}}
\put(1.3,0.8){\circle*{0.1}}

\put(3.2,0.2){\circle*{0.1}}
\put(3.8,0.2){\circle*{0.1}}
\put(3.2,0.8){\circle*{0.1}}
\put(3.8,0.8){\circle*{0.1}}

\put(1,0.5){\oval(2,1)}

\put(5.15,0.5){\vector(1,0){1}}

\put(8.2,0.2){\circle*{0.1}}
\put(8.8,0.2){\circle*{0.1}}
\put(8.2,0.8){\circle*{0.1}}
\put(8.8,0.8){\circle*{0.1}}

\put(10.7,0.2){\circle*{0.1}}
\put(11.3,0.2){\circle*{0.1}}
\put(10.7,0.8){\circle*{0.1}}
\put(11.3,0.8){\circle*{0.1}}

\put(8.5,0.5){\oval(2,1)}

\put(10.7,0.2){\line(1,0){0.6}}
\put(10.7,0.2){\line(0,1){0.6}}
\put(10.7,0.8){\line(1,0){0.6}}
\put(11.3,0.2){\line(0,1){0.6}}

\end{picture}

\caption{Constructing a component isomorphic to $C$.} \label{newfig1}
\end{figure}

For each graph $J \in \mathcal{J}_{n} \cap \mathcal{G}_{n}$,
delete all $|C|$ edges from a copy of $C$
(we have more than $r$ choices for this)
and insert them between $|C|$ isolated vertices
(we have at least $\left( ^{\alpha n}_{|C|} \right)$ choices for these)
to form a component isomorphic to $C$ (see Figure~\ref{newfig1}).
We can thus construct at least
$|\mathcal{J}_{n} \cap \mathcal{G}_{n}| r \left( ^{\alpha n}_{|C|} \right)$
graphs in $\mathcal{P}(n,m)$.

Given one of our constructed graphs,
there are at most $|C|+1$ components isomorphic to $C$,
since we have deliberately constructed one and we may have also created at most $|C|$ when we deleted the edges
(since each vertex in the deleted copy of $C$ might now be in a component isomorphic to $C$).
Thus, there are at most $|C|+1$ possibilities for which is the deliberately constructed component
and hence at most $|C|+1$ possibilities for which are the inserted edges.
There are then at most
$\left( ^{\phantom{i}n}_{|C|} \right) |C|! \leq n^{|C|}$
possibilities for where these edges were originally.
Thus, we have built each graph at most
$(|C|+1) n^{|C|}$ times
and so
\begin{eqnarray*}
|\mathcal{J}_{n} \cap \mathcal{G}_{n}|
& < & \frac{(|C|+1) n^{|C|} |\mathcal{P}(n,m)|}{ r \left( ^{\alpha n}_{\phantom{}|C|} \right)} \\
& = & (1+o(1)) \frac{(|C|+1)!}{r \alpha ^{|C|}} |\mathcal{P}(n,m)| \\
& < & \frac{\epsilon}{2} |\mathcal{P}(n,m)| \textrm{ for sufficiently large $n$,} \\
& & \textrm{ by definition of $r$.}
\end{eqnarray*}
Thus,
$|\mathcal{J}_{n}^{c}| \geq |\mathcal{J}_{n}^{c} \cap \mathcal{G}_{n}|
=|\mathcal{G}_{n}| - |\mathcal{J}_{n} \cap \mathcal{G}_{n}|
> \frac{\epsilon}{2} |\mathcal{P}(n,m)|$,
where $\mathcal{J}_{n}^{c}$ denotes $\mathcal{P}(n,m) \setminus \mathcal{J}_{n}$
(i.e.~the collection of graphs in $\mathcal{P}(n,m)$
without a set of more than $r$ edge-disjoint copies of $C$).

For $L \in \mathcal{J}_{n}^{c}$,
let $S(L)$ denote a maximal set of edge-disjoint copies of $C$ (so $|S(L)| \leq r$)
and, for $s \leq r$, let $\mathcal{J}_{n,s}$ denote the set of graphs in $\mathcal{P}(n,m)$ with $|S(L)|=s$.
Let us now consider the set $\mathcal{J}_{n,s}$.

\begin{figure} [ht]
\setlength{\unitlength}{1cm}
\begin{picture}(20,2.5)(0,0)

\put(0.5,2){\circle{1}}
\put(2.5,2){\circle{1}}
\put(4.5,2){\circle{1}}
\put(1.5,0.5){\circle{1}}
\put(3.5,0.5){\circle{1}}

\put(0.3,1.8){\line(1,0){0.4}}
\put(0.3,2.2){\line(1,0){0.4}}
\put(0.3,1.8){\line(0,1){0.4}}
\put(0.7,1.8){\line(0,1){0.4}}
\put(0.3,1.8){\circle*{0.1}}
\put(0.7,1.8){\circle*{0.1}}
\put(0.3,2.2){\circle*{0.1}}
\put(0.7,2.2){\circle*{0.1}}

\put(5.5,1.25){\vector(1,0){1}}

\put(7.5,2){\circle{1}}
\put(9.5,2){\circle{1}}
\put(11.5,2){\circle{1}}
\put(8.5,0.5){\circle{1}}
\put(10.5,0.5){\circle{1}}

\put(7.3,1.8){\circle*{0.1}}
\put(7.7,1.8){\circle*{0.1}}
\put(7.3,2.2){\circle*{0.1}}
\put(7.7,2.2){\circle*{0.1}}

\put(8,2){\circle*{0.1}}
\put(9,2){\circle*{0.1}}
\put(10,2){\circle*{0.1}}
\put(11,2){\circle*{0.1}}
\put(9.5,1.5){\circle*{0.1}}
\put(8.5,1){\circle*{0.1}}
\put(9,0.5){\circle*{0.1}}
\put(10,0.5){\circle*{0.1}}

\put(8,2){\line(1,0){1}}
\put(10,2){\line(1,0){1}}
\put(9,0.5){\line(1,0){1}}
\put(8.5,1){\line(2,1){1}}

\end{picture}

\caption{Destroying a copy of $C$.} \label{newfig2}
\end{figure}
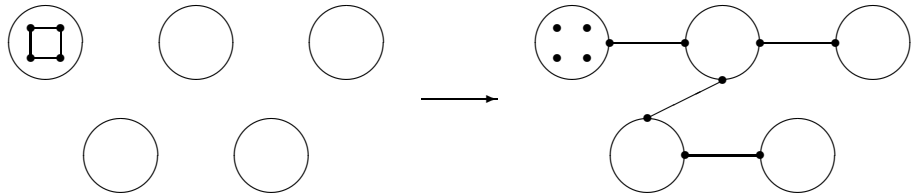

For each graph $L \in \mathcal{J}_{n,s}$,
delete all $|C|$ edges from a copy of $C$ that is in $S(L)$.
Note that the graphs will now have $|S|=s-1$,
by maximality of $S(L)$.
Clearly, we may insert an edge between any two vertices in different components
without introducing a copy of a cycle,
and, as in the proof of Lemma~\ref{add2},
we have at least $(1+o(1)) \frac{(1-A)^{2}}{2} n^{2}$ choices
for how to do this.
Thus, inserting $|C|$ edges,
we find that we may construct at least
$(1+o(1)) \frac{ \left( \frac{(1-A)^{2}}{2} \right)^{|C|} n^{2|C|} }{|C|!} |\mathcal{J}_{n,s}|$
graphs in $\mathcal{J}_{n,s-1}$ (see Figure~\ref{newfig2}).

Given one of our created graphs,
there are at most $\left( ^{\phantom{i}m}_{|C|} \right) \leq (An)^{|C|}$ possibilities for which edges were inserted
and at most
$\left( ^{\phantom{i}n}_{|C|} \right) |C|! < n^{|C|}$
possibilities for where they were originally.
Thus, we have built each graph at most $A^{|C|}n^{2|C|}$ times,
and so
$|\mathcal{J}_{n,s-1}| \geq
(1+o(1)) \frac{ \left( \frac{(1-A)^{2}}{2A} \right)^{|C|} }{|C|!} |\mathcal{J}_{n,s}|$.

For $s \leq r$,
let $\mathcal{J}_{n, \leq s}$ denote the set of graphs in $\mathcal{P}(n,m)$ with $|S| \leq s$.
Then, with $z = \frac{ \left( \frac{(1-A)^{2}}{2A} \right)^{|C|} }{|C|!}$, we have
$|\mathcal{J}_{n, \leq s-1}| \geq
(1+o(1)) \frac{z}{1+z} |\mathcal{J}_{n, \leq s}|$ for all $s \leq r$.
Thus, we obtain
\begin{eqnarray*}
|\mathcal{J}_{n, \leq 0}|
& \geq & (1+o(1))^{r} \left( \frac{z}{1+z} \right)^{r} |\mathcal{J}_{n, \leq r}| \\
& = & (1+o(1)) \left( \frac{z}{1+z} \right)^{r} |\mathcal{J}_{n}^{c}| \\
& \geq & (1+o(1)) \left( \frac{z}{1+z} \right)^{r} \frac{\epsilon}{2} |\mathcal{P}(n,m)|.
\end{eqnarray*}
But $\mathcal{J}_{n, \leq 0}$ is the set of graphs in $\mathcal{P}(n,m)$
without a copy of $C$,
and so we are done.~\phantom{qwerty}\begin{picture}(1,1)
\put(0,0){\line(1,0){1}}
\put(0,0){\line(0,1){1}}
\put(1,1){\line(-1,0){1}}
\put(1,1){\line(0,-1){1}}
\end{picture} \\
\\

In a moment,
we shall complete the unicyclic case by looking at what happens when $\frac{m}{n} \to 1$.
First, though,
we need to note two useful lemmas on $\kappa(P_{n,m})$,
the number of components of $P_{n,m}$:

\begin{Lemma}[\cite{dow}, Lemma 41] \label{cpt12}
Let $m(n) \in [n,3n-6]$ for all large $n$.
Then there exists a constant $c$ such that
\begin{displaymath}
\mathbf{P} \left[ \kappa (P_{n,m}) > \left \lceil \frac{cn}{\ln n} \right \rceil \right] = e^{- \Omega(n)}.
\end{displaymath}
\end{Lemma}
\textbf{Proof}
The proof follows that of Lemma 2.6 of~\cite{ger},
which deals with the case when $m=\lfloor qn \rfloor$ for fixed $q \in [1,3)$
(see~\cite{dow} for details).
\phantom{qwerty}
\begin{picture}(1,1)
\put(0,0){\line(1,0){1}}
\put(0,0){\line(0,1){1}}
\put(1,1){\line(-1,0){1}}
\put(1,1){\line(0,-1){1}}
\end{picture} \\

\begin{Lemma}[\cite{dow}, Proposition 50] \label{cpt31}
Let $q \in (0,1)$ be a constant and
let $m(n) \in [qn,n-1]$ for all large $n$.
Then there exists a constant $c=c(q)$ such that
\begin{displaymath}
\mathbf{P} \left[ \kappa (P_{n,m}) > \left \lceil \frac{cn}{\ln n} + n-m \right \rceil \right] = e^{- \Omega(n)}.
\end{displaymath}
\end{Lemma}
\textbf{Proof}
The proof follows that of Lemma 6.6 of~\cite{ger3},
which deals with the case when
$m = n - (\beta + o(1))(n / \ln n)$ for fixed $\beta > 0$
(see~\cite{dow} for details).
\phantom{qwerty}
\begin{picture}(1,1)
\put(0,0){\line(1,0){1}}
\put(0,0){\line(0,1){1}}
\put(1,1){\line(-1,0){1}}
\put(1,1){\line(0,-1){1}}
\end{picture} \\
\\

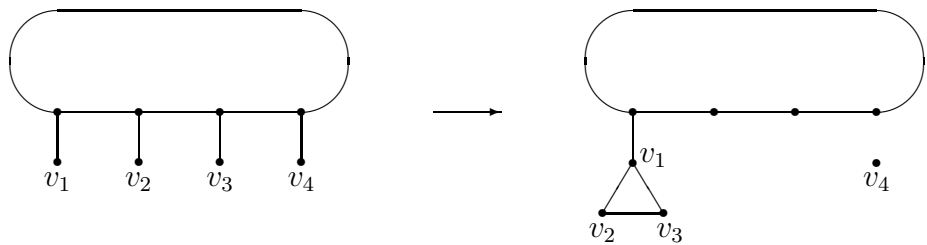
\begin{figure} [ht]
\setlength{\unitlength}{0.9cm}
\begin{picture}(20,3.35)(-0.75,-0.35)

\put(-0.05,0.75){\circle*{0.1}}
\put(1.15,0.75){\circle*{0.1}}
\put(2.35,0.75){\circle*{0.1}}
\put(3.55,0.75){\circle*{0.1}}

\put(1.75,2.25){\oval(5,1.5)}

\put(-0.05,0.75){\line(0,1){0.75}}
\put(-0.05,1.5){\circle*{0.1}}
\put(1.15,0.75){\line(0,1){0.75}}
\put(1.15,1.5){\circle*{0.1}}
\put(2.35,0.75){\line(0,1){0.75}}
\put(2.35,1.5){\circle*{0.1}}
\put(3.55,0.75){\line(0,1){0.75}}
\put(3.55,1.5){\circle*{0.1}}

\put(5.5,1.5){\vector(1,0){1}}

\put(10.25,2.25){\oval(5,1.5)}

\put(8.45,1.5){\circle*{0.1}}
\put(9.65,1.5){\circle*{0.1}}
\put(10.85,1.5){\circle*{0.1}}
\put(12.05,1.5){\circle*{0.1}}

\put(-0.25,0.4){$v_{1}$}
\put(0.95,0.4){$v_{2}$}
\put(2.15,0.4){$v_{3}$}
\put(3.35,0.4){$v_{4}$}

\put(12.05,0.75){\circle*{0.1}}
\put(8.45,0.75){\circle*{0.1}}
\put(8.45,0.75){\line(0,1){0.75}}

\put(8.45,0.75){\line(3,-5){0.45}}
\put(8.45,0.75){\line(-3,-5){0.45}}
\put(8,0){\circle*{0.1}}
\put(8.9,0){\circle*{0.1}}
\put(7.8,-0.35){$v_{2}$}
\put(8.8,-0.35){$v_{3}$}
\put(8.55,0.75){$v_{1}$}
\put(11.85,0.4){$v_{4}$}
\put(8,0){\line(1,0){0.9}}

\end{picture}

\caption{Constructing an appearance of $H$.} \label{subfig1}
\end{figure}

We may now prove our aforementioned result for when $\frac{m}{n} \to 1$:

\begin{Theorem} \label{unisub1}
Let $H$ be a (fixed) connected unicyclic graph
and let $m(n)=(1+o(1))n$.
Then
\begin{displaymath}
\mathbf{P}[P_{n,m} \textrm{ will have a copy of } H] \to 1
\textrm{ as } n \to \infty.
\end{displaymath}
\end{Theorem}
\textbf{Sketch of Proof}
We work with the set of graphs that don't have any copies of $H$.
By Lemma~\ref{pen5},
we may assume that we have lots of pendant edges,
so we may delete $|H|$ of these edges
and then use them with $|H|-1$ of the associated (now isolated) vertices
to convert another pendant edge into an appearance of $H$
(see Figure~\ref{subfig1}).
Note that we are also left with an extra isolated vertex.

By Lemmas~\ref{cpt12} and~\ref{cpt31},
we may assume that there are not very many isolated vertices.
Thus, since our original graphs had no appearances of $H$,
the amount of double-counting will be small,
and hence the size of our original set of graphs must have been small. \\
\\
\textbf{Full Proof}
Let $\mathcal{G}_{n}$ denote the set of graphs in $\mathcal{P}(n,m)$
with no copies of $H$,
and let $X$ denote the event that $P_{n,m} \in \mathcal{G}_{n}$.

Recall, from Lemma~\ref{pen5},
 that there exists an $\alpha >0$ such that
\begin{eqnarray}
\mathbf{P}[P_{n,m}
\textrm{ will have at least $\alpha n$ pendant edges}] \to 1
\textrm{ as } n \to \infty. \label{eq:u*}
\end{eqnarray}
Let $\mathcal{H}_{n}$ denote the set of graphs in $\mathcal{P}(n,m)$
with at least $\alpha n$ pendant edges,
and let $Y$ denote the event that $P_{n,m} \in \mathcal{H}_{n}$.

Also, note that by Lemmas~\ref{cpt12} and~\ref{cpt31},
there exists a $c>0$ such that
\begin{eqnarray}
\mathbf{P} \left[ \kappa \left( P_{n,m} \right)
> 2 \max \left \{ \frac{cn}{\ln n}, n-m \right \} \right] = e^{- \Omega(n)}. \label{eq:udag}
\end{eqnarray}
Let $x=x(n)= 2 \max \{ \frac{cn}{\ln n}, n-m \} = o(n)$,
let $\mathcal{I}_{n}$ denote the number of graphs in $\mathcal{P}(n,m)$ with at most $x$ components,
and let $Z$ denote the event that $P_{n,m} \in \mathcal{I}_{n}$.

Note that
$\mathbf{P}[X] \leq \mathbf{P}[X \cap Y \cap Z] + \mathbf{P} \left[\bar{Y}\right] + \mathbf{P}\left[\bar{Z}\right]
\to \mathbf{P}[X \cap Y \cap Z] \textrm{ as } n \to \infty$,
by (\ref{eq:u*}) and (\ref{eq:udag}).
Thus, it suffices to show that $\mathbf{P}[X \cap Y \cap Z] \to 0$ as $n \to \infty$.

Given a graph in $\mathcal{G}_{n} \cap \mathcal{H}_{n} \cap \mathcal{I}_{n}$,
let us choose $|H|+1$ pendant edges
(we have at least $\left( ^{\phantom{w} \alpha n}_{|H|+1} \right)$ choices for these).
We shall use these edges to create an appearance of $H$.

Out of our chosen $|H|+1$ pendant edges,
let the edge incident with the lowest labelled vertex of degree $1$ be called `special',
and let this associated lowest labelled vertex of degree $1$ be called the `root'.
Let us delete our $|H|$ non-special pendant edges to create at least $|H|$ isolated vertices.

We may choose $|H|-1$ of these newly isolated vertices in such a way
that no two were adjacent in the original graph
(i.e.~we don't choose two vertices from the same pendant edge,
even if that is possible).
We may then use these chosen isolated vertices,
together with the root,
to construct an appearance of $H$
by inserting $|H|$ edges appropriately
(see Figure~\ref{subfig1}).

We shall now consider the amount of double-counting:

Recall that our original graphs contained no copies of $H$
and note that we cannot have created any copies by deleting edges.
As observed in the proof of Theorem 4.1 of~\cite{mcd},
an appearance of $H$ meets at most $|H|-1$ other appearances of $H$
(since there are at most $|H|-1$ cut-edges in $H$
and each of these can have at most one `orientation' that provides an appearance of $H$).
Thus,
when we deliberately constructed our appearance of $H$
we can have only created at most $|H|$ appearances of $H$ in total.

Therefore, given one of our constructed graphs,
there are at most $|H|$
possibilities for which is our deliberately created appearance of $H$.
We may then recover the original graph by deleting the $|H|$ edges from this appearance,
joining the $|H|-1$ non-root vertices back to the rest of the graph
(at most $n^{|H|-1}$ possibilities),
and joining the correct isolated vertex back to the rest of the graph
(at most $in$ possibilities,
where $i$ denotes the number of isolated vertices in the constructed graph).

We know that the number of components in the original graph was at most $x$,
and each time we deleted an edge we can have only increased the number of components by at most $1$.
Thus, the number of components in our constructed graph is at most $x+|H|$,
and so $i \leq x+|H|$.

Thus, we have built each graph at most
$|H|n^{|H|-1}(x+|H|)n$ times.

Therefore,
\begin{eqnarray*}
\mathbf{P}[X \cap Y \cap Z]
& = & \frac{|\mathcal{G}_{n} \cap \mathcal{H}_{n} \cap \mathcal{I}_{n}|}{|\mathcal{P}(n,m)|} \\
& \leq & \frac{|H|n^{|H|}(x+|H|)}
{\left( ^{\phantom{w} \alpha n}_{|H|+1} \right) } \\
& = & \frac{o \left( n^{|H|+1} \right) }{\Theta \left( n^{|H|+1} \right) }, \textrm{ since } x=o(n)\\
& \to & 0 \textrm{ as } n \to \infty.
\phantom{qwerty}\begin{picture}(1,1)
\put(0,0){\line(1,0){1}}
\put(0,0){\line(0,1){1}}
\put(1,1){\line(-1,0){1}}
\put(1,1){\line(0,-1){1}}
\end{picture}
\end{eqnarray*}
\\

\begin{figure} [ht]
\setlength{\unitlength}{0.97cm}
\begin{picture}(20,1.5)(-0.25,0)

\put(1.8,0.4){\line(1,0){0.9}}
\put(1.8,0.4){\line(3,2){0.45}}
\put(1.8,0.4){\line(3,5){0.45}}
\put(2.25,0.7){\line(0,1){0.45}}
\put(2.7,0.4){\line(-3,2){0.45}}
\put(2.7,0.4){\line(-3,5){0.45}}
\put(1.8,0.4){\circle*{0.1}}
\put(2.7,0.4){\circle*{0.1}}
\put(2.25,0.7){\circle*{0.1}}
\put(2.25,1.15){\circle*{0.1}}

\put(2.25,0.75){\oval(5,1.5)}

\put(5.5,0.75){\vector(1,0){1}}

\put(8.25,0.4){\line(1,0){0.5}}
\put(8.25,0.75){\line(1,0){0.5}}
\put(8.25,1.1){\line(1,0){0.5}}
\put(10.75,0.4){\line(1,0){0.5}}
\put(10.75,0.75){\line(1,0){0.5}}
\put(10.75,1.1){\line(1,0){0.5}}

\put(9.3,0.4){\circle*{0.1}}
\put(10.2,0.4){\circle*{0.1}}
\put(9.75,0.7){\circle*{0.1}}
\put(9.75,1.15){\circle*{0.1}}

\put(9.75,0.75){\oval(5,1.5)}

\end{picture}

\caption{Redistributing the edges of $H$.} \label{msubfig}
\end{figure}
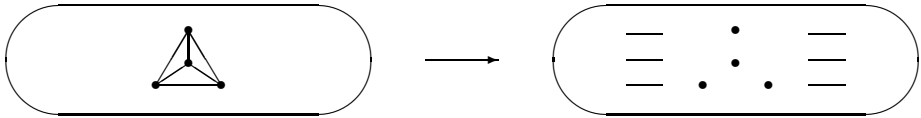

It now only remains to look at the case when $\frac{m}{n} \leq 1+o(1)$
and $H$ is a connected multicyclic graph.
If $\limsup \frac{m}{n} < 1$,
we may prove the following:

\begin{Theorem} \label{msub3}
Let $H$ be a (fixed) multicyclic connected planar graph
and let $m(n)$ satisfy
$\limsup_{n \to \infty} \frac{m}{n} < 1$.
Then
\begin{displaymath}
\mathbf{P}[P_{n,m} \textrm{ will have a copy of }H] \to 0 \textrm{ as } n \to \infty.
\end{displaymath}
\end{Theorem}
\textbf{Proof}
Let $\mathcal{G}_{n}$ denote the set of graphs in $\mathcal{P}(n,m)$ with a copy of $H$.
For each graph $G \in \mathcal{G}_{n}$, let us delete all $e(H)$ edges from a copy of $H$
and then insert these edges back into the graph (see Figure~\ref{msubfig}).
By Lemma~\ref{add2},
we have $\Theta \left( n^{2e(H)} \right)$ ways to do this, maintaining planarity.

There are then
$\left(\left(^{\phantom{i}n}_{|H|}\right)\frac{|H|!}{|\textrm{\scriptsize{Aut}}(H)|} \right) = O\left( n^{|H|} \right)$
possibilities
for where the copy of $H$ was originally
and $\left(^{\phantom{n}m} _{e(H)} \right) = O\left( n^{e(H)} \right)$
possibilities for which edges were inserted,
so we have built each graph $O\left( n^{(|H|+e(H))} \right)$ times.

Therefore,
$\frac{|\mathcal{G}_{n}|}{|\mathcal{P}(n,m)|} = O \left( \frac{n^{(|H|+e(H))}}{n^{2e(H)}} \right)
\to 0$,
since $e(H)>|H|.$
\phantom{qwerty}
\setlength{\unitlength}{.25cm}
\begin{picture}(1,1)
\put(0,0){\line(1,0){1}}
\put(0,0){\line(0,1){1}}
\put(1,1){\line(-1,0){1}}
\put(1,1){\line(0,-1){1}}
\end{picture}
\\
\\

By using more precise estimates of add$(n,m)$,
it is possible to slightly increase the upper bound imposed on $m$ in Theorem~\ref{msub3}
to include some functions $m(n)$ for which $\frac{m}{n}$ converges to $1$ slowly
(see Theorems 68 and 70 of~\cite{dow}).
However,
we do not yet have a result that covers the whole region $m=(1+o(1))n$. \\

\section*{Acknowledgements}

I am very grateful to Colin McDiarmid for his advice and helpful suggestions. \\

\end{document}